\documentclass[11pt,reqno]{amsart}
\usepackage{fullpage,amsfonts,amsmath,amscd,amssymb}

\theoremstyle{plain}
\newtheorem*{trm}{Theorem}
\newtheorem*{lem}{Lemma}
\newtheorem*{prop}{Proposition}
\newtheorem*{thm}{Theorem}

\theoremstyle{remark}

\newtheorem*{rem}{Remark}

\newcommand{\C}{\mathbb{C}}
\newcommand{\spec}{\text{spec}}
\newcommand{\Cas}{\text{Cas}}
\newcommand{\mat}{\text{Mat}}

\newcommand{\EO}{\mathcal{O}_{\epsilon}[G]}
\newcommand{\CO}{\mathcal{O}[G]}
\title[Poisson orders]{Poisson orders, symplectic reflection algebras and
representation theory}
\author{Kenneth A.Brown}
\author{Iain Gordon}
\thanks{The second author thanks the Nuffield Foundation for
the financial support provided by grant number NAL/00625/G}
\email{kab@maths.gla.ac.uk, ig@maths.gla.ac.uk}
\begin{document}

\begin{abstract} We introduce a new class of algebras called Poisson
orders. This class includes the symplectic reflection algebras of
Etingof and Ginzburg, many quantum groups at roots of unity, and
enveloping algebras of restricted Lie algebras in positive
characteristic. Quite generally, we study this class of algebras
from the point of view of Poisson geometry, exhibiting connections
between their representation theory and some well-known geometric
constructions. As an application, we employ our results in the
study of symplectic reflection algebras, completing work of
Etingof and Ginzburg on when these algebras are finite over their
centres, and providing a framework for the study of their
representation theory in the latter case.
\end{abstract}\maketitle
\section{Introduction}
\subsection{}\label{start} It is a truth universally acknowledged, that a geometric
viewpoint benefits representation theory. Evidence is found in the
success of the orbit method and the theory of characteristic
varieties, whose tool is symplectic geometry, and the theory of
algebras which are finite as modules over their centres, in which
algebraic geometry is pervasive.

In this paper we introduce a new class of algebras,
\textit{Poisson orders}, which exhibit key features from both of
the above exemplars. Poisson orders include as special cases the
$t=0$ case of symplectic reflection algebras (whose definition we
recall in (\ref{sradefs})), enveloping algebras of restricted Lie
algebras in positive characteristic, and quantised enveloping
algebras and function algebras at roots of unity. We prove several
basic theorems which relate the representation theory and
algebraic structure of Poisson orders to symplectic geometry and,
as an extended example, study symplectic reflection algebras.
\subsection{}
A Poisson order is an affine $\mathbb{C}$-algebra $A$, finitely
generated as a module over a central subalgebra $Z_0$, together
with the datum of a $\mathbb{C}$-linear map $D:Z_0 \longrightarrow
\text{Der}_\mathbb{C}(A)$, which satisfies the Leibniz identity
and makes $Z_0$ a Poisson algebra by restriction. The details of
the definition are given in (\ref{setup}). An elementary but
important example is the skew group algebra $A = \mathcal{O}(V)
\ast \Gamma$, where $\mathcal{O}(V)$ is the algebra of polynomial
functions on a symplectic $\mathbb{C}$-vector space $V$ and
$\Gamma$ is a finite subgroup of $\textrm{Sp}(V)$. Then $A$ is a
Poisson order with $Z_0 = \mathcal{O}(V)^{\Gamma}$, the centre of
$A$. The Poisson structure is induced from the symplectic form on
$V$.

Another particularly rich source of examples is quantisation,
recalled in (\ref{alg1}). This yields in particular the examples
of quantum groups and Lie algebras mentioned in (\ref{start}), but
it is also the mechanism whereby symplectic reflection algebras in
the $t = 0$ case are given the structure of Poisson orders,
(\ref{sraPois}).

\subsection{}\label{stratover}
Let $Z_0 \subseteq A$ be a Poisson order and set $\mathcal{Z} =
\text{Maxspec}Z_0$, a Poisson variety. The Poisson structure
induces several stratifications on $\mathcal{Z}$, as follows.

\begin{itemize}
\item The \textit{rank stratification}, where $\mathcal{Z}$
is stratified by the rank of the Poisson bracket at a point.
\item
The stratification by \textit{symplectic cores}, refining the rank
stratification. We expect this stratification is the finest
possible algebraic stratification in which the Hamiltonian vector
fields, $D(z)$ for $z\in Z_0$, are tangent.
\item
The stratification (when $A$ is a $\C$-algebra) by
\textit{symplectic leaves}, a further refinement, which
potentially leaves the algebraic category. This is a differential
geometric notion, which produces the well-known foliation by
symplectic leaves when $\mathcal{Z}$ is smooth.
\end{itemize}
We introduce and discuss the relations between these
stratifications in Section 3. When the stratification by
symplectic cores and by symplectic leaves agree, we call the
Poisson bracket on $Z_0$ \textit{algebraic}. In favourable
circumstances, for instance if there is only a finite number of
symplectic leaves, these three stratifications essentially agree.
In general, however, we know of no procedure to deduce when a
bracket is algebraic.
\subsection{}
We show in this paper that the symplectic cores have a profound
influence on the representation theory of $A$. Indeed, the
surjective map
\[
\mathcal{I}rr(A) \longrightarrow \mathcal{Z},
\]
sending an irreducible $A$-module, $S$, to its $Z_0$-character
$\text{Ann}_{Z_0}(S)$, \textit{has constant fibres over the
symplectic cores of $\mathcal{Z}$}. In fact more is true. Given $x
\in \mathcal{Z}$, let $A_x$ be the finite dimensional algebra
obtained by taking the quotient of $A$ by the ideal
$\mathfrak{m}_xA$. We show in Theorem \ref{algIM} that {\em if
$x,y \in\mathcal{Z}$ and $x$ and $y$ belong to the same core, the
algebras $A_x$ and $A_y$ are isomorphic}. Moreover, as we explain
in Section 5, this isomorphism can be globalised, yielding a sheaf
of \textit{Azumaya} algebras over each core.

Results of this sort were first proved in the case of quantised
function algebras at a root of unity in \cite{DeCL} and
\cite{DeCP}. The more general result quoted above can also be
fruitfully applied to symplectic reflection algebras. As we show
in Section 7, {\em these algebras are sufficiently well-behaved to
ensure that the stratifications in (\ref{stratover}) agree}. In
Section 7 we also take the opportunity to complete work begun in
\cite{EG} {\em by determining when a symplectic reflection algebra
is a finite module over its centre, and showing that the centre is
$\C$ in all other cases}. These results show that it is a
fundamental problem to determine the symplectic leaves of the
centres of symplectic algebras. We give a complete answer to this
in the ``degenerate" case of the skew group algebra
$\mathcal{O}(V) \ast \Gamma$, so answering a question of
\cite{AF}, and we provide some preliminary information in the
general case, in (\ref{0case})-(\ref{gencase}).
\subsection{}
The paper is organised as follows. In Section \ref{pss} we
introduce the basic definition of a Poisson order and its Poisson
ideals. We follow this in Section \ref{stt} by studying the
various stratifications of $\mathcal{Z}$, pausing to discuss
relations to the \textit{Dixmier-Moeglin equivalence}. In Sections
\ref{iss} and \ref{az} we prove our results relating the
representation theory of $A$ with the symplectic cores of
$\mathcal{Z}$. Section \ref{qfa} briefly recalls the (motivating)
case of quantised function algebras at roots of unity. Finally, in
Section \ref{sra} we discuss symplectic reflection algebras.

\subsection{}
Although we work over a field of characteristic zero, which in
later sections is assumed to be $\mathbb{C}$, most of the
algebraic results can be generalised to arbitrary algebraically
closed fields. However, there are serious obstacles to be overcome
in the attempt to generalise the geometric picture to other
(especially positive characteristic) fields. Given that the
enveloping algebra of a restricted Lie algebra comprises part of
our motivation, this is an important issue requiring attention.

\section{Poisson orders}\label{pss}

\subsection{Definition and notation}
\label{setup} Let $k$ be a field of characteristic $0$ and let $A$
be an affine $k$-algebra, finitely generated as a module over a
central subalgebra $Z_0$. By the Artin-Tate Lemma
\cite[13.9.10]{McC-Rob}, $Z_0$ is an affine $k$-algebra, whose
(Krull) dimension we shall fix as $d$ throughout. Suppose that
there is a linear map
\[
D: Z_0 \longrightarrow \text{Der}_k(A): z \mapsto D_z,
\]
 satisfying
\begin{enumerate}
\item $Z_0$ is stable under $D(Z_0)$;
\item the resulting bracket
$ \{ - , -\} : Z_0\times Z_0\longrightarrow Z_0, $ defined by $\{
z, z' \}= D_z(z')$, imposes a structure of Poisson $k$-algebra on
$Z_0$ - that is, $D_{zz'} = zD_{z'} + z'D_{z}$ as derivations of
$Z_0$, and $(Z_0,\{ - , -\})$ is a $k$-Lie algebra.
\end{enumerate}
Then we shall say that $A$ is a \textit{Poisson $Z_0$-order}. We
can also express (2) as saying that $\mathcal{Z}: =
\text{Maxspec}(Z_0)$ is a \textit{Poisson variety} with respect to
the bracket defined in (2). The \textit{algebra of Casimirs} of a
Poisson algebra $Z_0$ is $\Cas (Z_0) := \{ z \in Z_0 : \{z,Z_0 \}
= 0 \}.$

Throughout this section we'll  assume that $A$ is a Poisson
$Z_0$-order with associated linear map $D$.

\subsection{Quantisation}\label{alg1}
Here is one important mechanism giving rise to a Poisson order.
 Let $\hat{A}$ be a $k$-algebra, $\hat{Z}$ a
subalgebra, and $t\in \hat{Z}$ a central non-zero divisor. Assume
that $Z_0=\hat{Z}/t\hat{Z}$ is an affine central subalgebra of
$A=\hat{A}/t\hat{A}$, and that $A$ is a finitely generated
$Z_0$-module. Let $\pi:\hat{A}\longrightarrow A$ be the quotient
map. Fix a $k$-basis $\{z_i : i \in \mathcal{I} \}$ of $Z_0$, and
lift these elements of $Z_0$ to elements $\hat{z}_i$ of $\hat{Z}$.

Given $i \in \mathcal{I}$ there is a derivation of $A$, denoted
$D_{z_i}$, defined by
\[ D_{z_i} (a) = \pi([\hat{z}_i, \hat{a}]/t), \]
where $\hat{a}\in \hat{A}$ is a preimage under $\pi$ of $a$. Now
define $D:Z_0 \longrightarrow \text{Der}(A)$ by extending the
above definition $k$-linearly. This construction satisfies the
hypotheses of (\ref{setup}), \cite{H}. Observe that alternative
choices (for example of liftings of $\{z_i \}$) would yield the
same outcomes, to within inner derivations of $A$.

\subsection{Filtered and graded algebras}\label{alg2}
An important variant of the above is the following. Let $A$ be an
$\mathbb{N}$-filtered algebra whose $i^{\text{th}}$-filtered piece
is denoted $F^iA$. Let $Z$ be a subalgebra of $A$, and give it the
induced filtration. Denote the associated graded rings of $Z$ and
$A$ by $\text{gr}Z$ and $\text{gr}A$ respectively. Suppose that
$\text{gr}Z$ is an affine central subalgebra of $\text{gr}A$, such
that $\text{gr}A$ is a finitely generated $\text{gr}Z$-module. Let
$\sigma_i:F^iA \longrightarrow \text{gr}A$ be the
$i^{\text{th}}$-principal symbol map, sending an element of
$F^iA\setminus F^{i-1}A$ to its leading term. Given a homogeneous
element of $\text{gr}Z$, say $\sigma_m(z)$, there is a
well-defined derivation of $\text{gr}A$, denoted
$D_{\sigma_m(z)}$, given, for a homogeneous element $\sigma_n (a)$
of $\text{gr}(A)$, by \begin{equation} \label{p1}
D_{\sigma_m(z)}(\sigma_n(a)) = \sigma_{m+n-1}([z,a]).
\end{equation} Extending this linearly yields a mapping
$D:\text{gr}Z \longrightarrow \text{Der}(\text{gr}A)$, satisfying
the hypotheses of (\ref{setup}). Naturally, we'll call a Poisson
bracket on a commutative graded algebra $H$ {\em homogeneous} if
the bracket of any two homogeneous elements of $H$ is also
homogeneous; in this case, if whenever $h \in H_i$ and $g \in H_j$
we have $\{h,g\} \in H_{i+j+d}$, we shall say that the bracket has
{\em degree} $d$. Thus the bracket defined in (\ref{p1}) is
homogeneous of degree $-1$.

To see that this definition is really a special case of
(\ref{alg1}), form the Rees algebras $\hat{A} = \oplus_i F^iAt^i
\subseteq A[t]$ and $\hat{Z} = \oplus_i F^iZt^i\subseteq Z[t]$,
where $t$ is a central non-zero divisor. It can easily be checked
that we recover the derivations in (\ref{p1}) from the
construction in (\ref{alg1}).

\subsection{Poisson ideals and subsets}
\label{pcl} A two-sided ideal $I$ of the Poisson $Z_0$-order $A$
(respectively $J$ of $Z_0$) is called \textit{Poisson} if it is
stable under $D(Z_0)$. Thanks to \cite[3.3.2]{Dix} if $I$
(respectively $J$) is Poisson then so too are both $\sqrt{I}$
(respectively $\sqrt{J}$) and the minimal prime ideals of $A$
(respectively $Z_0$) over $I$ (respectively $J$). We shall denote
the space of prime Poisson ideals of $Z_0$, with the topology
induced from the Zariski topology on $\spec( Z_0)$, by
$\mathcal{P}-\spec (Z_0)$. Clearly, if $I$ is a Poisson ideal of
$Z_0$ then there is an induced structure of Poisson algebra on
$Z_0/I$.  For a semiprime ideal $I$ of $Z_0$ (respectively of $A$)
we write $\mathcal{V}(I)$ for the closed subset of $\mathcal{Z}$
(respectively of $\text{Maxspec}(A)$) defined by $I$. A closed
subset $\mathcal{V}(I)$ of $\mathcal{Z}$ is \textit{Poisson
closed} if its defining ideal is Poisson.

Extension ($e:J\longrightarrow JA$) and contraction ($c:I \longrightarrow I\cap Z_0$)
are mappings between the ideals of $Z_0$ and $A$, which map Poisson ideals to Poisson ideals.
It is an easy exercise using
Going Up \cite[10.2.10(ii)]{McC-Rob} to show that
\begin{equation*}
\label{g1} c \circ e \text{ is the identity on semiprime ideals of } Z_0.
\end{equation*}
Clearly, therefore, if $I$ is a semiprime Poisson ideal of $Z_0$
then there is an induced structure of Poisson $Z_0/I$-order on
$A/IA,$ and if $J$ is a Poisson ideal of $A$ then there is an
induced structure of Poisson $Z_0/J \cap Z_0$-order on $A/J.$

\section{Stratifications}\label{stt}
\subsection{The rank stratification}
\label{rank} Throughout Section 3 we consider an affine
commutative Poisson $k$-algebra $Z_0$. Let $\{ z_1,\ldots ,z_m\}$
be a generating set for $Z_0$ and define the $m\times m$-skew
symmetric matrix $\mathcal{M} = \big(\{ z_i,z_j\}\big ) \in
M_m(Z_0)$.
 The \textit{rank} of the
Poisson structure at $\mathfrak{m} \in \mathcal{Z}$, denoted
$\text{rk}(\mathfrak{m})$, is defined to be the rank of the matrix
$\mathcal{M}(\mathfrak{m}) = \big(\{ z_i,z_j\} + \mathfrak{m} \big
) \in M_m(k)$. It is independent of the choice of generators,
\cite[2.6]{van}. For each non-negative integer $j$ we define
$$\mathcal{Z}^o_j(Z_0) :=  \{ \mathfrak{m} \in  \mathcal{Z} : \text{rk}(\mathfrak{m}) = j\}$$
and
$$\mathcal{Z}_j(Z_0) := \{ \mathfrak{m} \in
\mathcal{Z} : \text{rk}(\mathfrak{m}) \leq j\}.$$ We'll simply
write $\mathcal{Z}^o_j$ and $\mathcal{Z}_j$ when the algebra
involved is clear from the context.

\begin{lem} Retain the above notation.
\begin{enumerate}
\item $\mathcal{Z}_j$ is a closed subset of $\mathcal{Z}$, with
$$ \mathcal{Z}_0 \subseteq \mathcal{Z}_1 \subseteq \ldots \subseteq \mathcal{Z}_d =
\mathcal{Z}.$$
\item $\mathcal{Z}_j$ is a Poisson subset of $\mathcal{Z}$, at least
when $k = \C.$
\item $\mathcal{Z}_j = \cup_{i \leq j} \mathcal{Z}^o_i =
\mathcal{Z}_{j_0}$, where $j_0 := \max \{ i : i \leq j,
\mathcal{Z}^o_i \neq \emptyset \}.$
\item The sets $\mathcal{Z}^o_i$ are locally closed: if $\mathcal{Z}^o_i$
is non-empty then $\overline{\mathcal{Z}^o_i}$ is equal to the
union of some of the irreducible components of $\mathcal{Z}_i$.
\item If $\mathcal{Z}_i^o\neq \emptyset$ then $\dim \mathcal{A} \geq i$
 for each irreducible component $\mathcal{A}$ of
$\mathcal{Z}^o_i.$
\item Suppose that $k = \C$. If the irreducible component $\mathcal{A}$ of $\mathcal{Z}^o_i$
has $\dim \mathcal{A} = i$ then $\mathcal{A}$ is contained in the
smooth locus of its closure $\overline{\mathcal{A}}$.
\end{enumerate} \end{lem}

\begin{proof} (1) It's clear that $\mathcal{Z}_j$ is closed, since $\mathcal{Z}_j$ is defined by
the vanishing of all minors of $\mathcal{M}$ of order $j + 1$.
That $\mathcal{Z}_i \subseteq \mathcal{Z}_{i+1}$ is clear from the
definition; and that $\mathcal{Z}_d = \mathcal{Z}$ is a
consequence of (5), which we prove below.

(2) Assume that $k = \C$. To show that $\mathcal{Z}_j$ is Poisson
we argue by noetherian induction, and so assume that the result is
true for all proper Poisson factors of $Z_0$. Let $J$ be the
defining ideal of $\mathcal{Z}_j$, and let $\mathcal{Z}_j$ have
irreducible components $\mathcal{A}_1 = \mathcal{V}(P_1), \ldots
,\mathcal{A}_t = \mathcal{V}(P_t)$. Thus $J$ is Poisson if and
only if each $P_i$ is Poisson, by (\ref{pcl}). We thus aim to show
that $P := P_1$ is Poisson. By (\ref{pcl}) again, we can assume
that $Z_0$ is a domain. Let $I$ be the defining ideal of the
singular locus of $\mathcal{Z}$, (so $I = Z_0$ if $Z_0$ is
smooth). Let $0 \neq y \in I.$ Then $Z_0[y^{-1}]$ is a Poisson
algebra and if $P[y^{-1}]$ is a proper ideal of $Z_0[y^{-1}]$ then
it defines a component of $\mathcal{Z}_j(Z_0[y^{-1}]).$ Since
$Z_0[y^{-1}]$ is smooth, $P[y^{-1}]$ is Poisson by \cite[Corollary
2.3]{Pol}.

Define $\hat{P} := \bigcap_{0 \neq y \in I} P[y^{-1}] \cap Z_0$,
so that $\hat{P}$ is a Poisson ideal of $Z_0$ containing $P$, with
$\hat{P}I \subseteq P.$ Since $P$ is prime either (i) $\hat{P} =
P$ and so $P$ is Poisson, or (ii) $I \subseteq P.$ Now $I$ is a
Poisson ideal by \cite[Corollary 2.4]{Pol}, and $I$ is non-zero by
definition. So if (ii) holds we can pass to $Z_0/I$ and invoke our
induction hypothesis. This completes the proof in either case.

(3), (4) If $\mathcal{Z}^o_i$ is non-empty then it is precisely
the subset of $\mathcal{Z}_i$ on which some minor of $\mathcal{M}$
of order $i$ is not identically zero.

(5) Suppose that $\mathcal{Z}^o_i$ is non-empty, and let
$\mathcal{A}$ be one of its
 irreducible components. Since $\mathcal{A}$
is open in its closure $\mathcal{T},$  $\mathcal{A}$ contains a
point $\mathfrak{m}$ which is smooth in $\mathcal{T}$. Let $J$ be
the defining ideal of $\mathcal{T}$, a Poisson ideal by (2) and
(\ref{pcl}). Note that $\mathcal{M}$ defines a bilinear form of
rank $i$ on $\mathfrak{m}/\mathfrak{m}^2$, and since $J$ is
Poisson, $J + \mathfrak{m}^2/\mathfrak{m}^2$
 is contained in the radical of this form. Therefore,
$$ i \leq  \dim(\mathfrak{m}/J + \mathfrak{m}^2) = \dim(\mathcal{T}) =
\dim(\mathcal{A}).$$

(6) Assume now that $\dim \mathcal{A} = i$, and suppose that
$\mathcal{A}$ contains a singular point of
$\overline{\mathcal{A}}.$ Then the intersection of $\mathcal{A}$
with the singular locus of $\overline{\mathcal{A}}$ is a non-empty
open subset of the singular locus, and so contains a smooth point
$\mathfrak{n}$ of that singular locus. Let $K$ be the ideal of
$Z_0$ defining the singular locus of $\overline{\mathcal{A}}$ , so
$K$ is a Poisson ideal by (2) and \cite[Corollary 2.4]{Pol}.
Hence, the image of $K$ in $\mathfrak{n}/\mathfrak{n}^2$ is in the
radical of the form induced on $\mathfrak{n}/\mathfrak{n}^2$. Thus
$$ i = \dim \overline{\mathcal{A}} > \text{Krull dim}(Z_0/K) =
\dim_{\C}(\mathfrak{n}/\mathfrak{n}^2 + K) \geq i,$$ a
contradiction.
 \end{proof}

In the light of the lemma we define the \textit{rank
stratification} of $\mathcal{Z}$ to be the disjoint union of
locally closed subsets
\begin{equation}
\label{strat1} \mathcal{Z}= \coprod_{j} \mathcal{Z}^o_{j}.
\end{equation}

\subsection{Dixmier-Moeglin equivalence for Poisson ideals}
\label{Dix-Moeg} Continue with the notation and assumptions
already introduced for $Z_0$. Given an ideal $I$ of $Z_0$ we
define the \textit{Poisson core} of $I$ to be the biggest Poisson
ideal of $Z_0$ contained in $I$. This exists since the sum of two
Poisson ideals is Poisson. We denote the Poisson core of $I$ by
$\mathcal{P}_{Z_0}(I)$, or just by $\mathcal{P}(I)$ when the
algebra is clear from the context.   If $I$ is prime then by
(\ref{pcl}) so also is $\mathcal{P}(I)$. A prime ideal
$\mathfrak{p}$ of $Z_0$ is \textit{Poisson primitive} if
$\mathfrak{p} = \mathcal{P}(\mathfrak{m})$ for a maximal ideal
$\mathfrak{m}$ of $Z_0$. Since every prime ideal of $Z_0$ is an
intersection of maximal ideals,
\begin{equation}\label{P-nullstell} \textit{every prime Poisson ideal of $Z_0$
is an intersection of Poisson primitive ideals.}
\end{equation}
By analogy with well-known results of Dixmier and Moeglin for
certain noncommutative affine algebras such as enveloping algebras
of finite dimensional Lie algebras \cite[Chapter 9]{McC-Rob} we
are led to consider the following hypotheses which can be imposed
on a prime Poisson ideal $\mathfrak{p}$ of $Z_0$:
\renewcommand{\theenumi}{\Alph{enumi}}
\begin{enumerate}
\item  $\mathfrak{p}$ is Poisson primitive;
\item  $\mathfrak{p}$ is locally closed in $\mathcal{P}-\spec(Z_0)$;
\item  $\Cas (Z_0/\mathfrak{p})$ is an algebraic extension of $k$;
\item  $\Cas\big(\text{Fract}(Z_0/\mathfrak{p})\big)$ is an algebraic extension of $k$.
\end{enumerate}
\renewcommand{\theenumi}{\arabic{enumi}}
The following lemma is due to Oh \cite[1.7(i),1.10]{Oh}.

\begin{lem} Let $\mathfrak{p} \in \mathcal{P}-\spec(Z_0)$, and fix hypotheses and notation as
above. Then $(B) \Longrightarrow (A) \Longrightarrow (D)
\Longrightarrow (C)$.
\end{lem}

\noindent \textbf{Remark:} In general (C) does not imply (A) or
(B). For consider $Z_0 = \C [x,y,z]$, with $\{x,z\} = x, \{y,z\} =
y, \{x,y\} = 0.$ If $F$ is a Casimir, then $\partial F /
\partial z = x \partial F /\partial x + y\partial F /\partial y =
0,$ so it follows easily that $F$ is a scalar. But one also
calculates that $\mathcal{P}-\spec (Z_0)$ consists of $\{ 0 \}$,
$\langle x - ay \rangle,$ for $a \in \C,$ and $\langle y \rangle$
and $\langle x,y \rangle$. All of these are thus Poisson primitive
except for $\{ 0 \}$, which also fails to be locally closed in
$\mathcal{P}-\spec (Z_0)$.

Note, however, that in the above example
$\Cas\big(\text{Fract}(Z_0)\big) = \C (xy^{-1}).$ This motivates
the following question, to which we will return in (\ref{D-M2})
and (\ref{algIM}). \smallskip

\noindent \textbf{Question:} Are properties (A), (B) and (D)
always equivalent for an affine Poisson $\C$-algebra $Z_0$?

\subsection{The stratification by symplectic cores}
\label{cores}
 We again
focus on $Z_0$ in this paragraph, and continue with the hypotheses
as above. We define a relation $\sim$ on $\mathcal{Z}$ by:
$$ \mathfrak{m} \sim \mathfrak{n} \Longleftrightarrow
\mathcal{P}(\mathfrak{m})=\mathcal{P}(\mathfrak{n}). $$ Clearly
$\sim$ is an equivalence relation; we denote the equivalence class
of $\mathfrak{m}$ by $\mathcal{C}(\mathfrak{m})$, so that
\begin{equation}
\label{strat2} \mathcal{Z} \quad = \quad \coprod
\mathcal{C}(\mathfrak{m}).
\end{equation}
The set $\mathcal{C}(\mathfrak{m})$ is called the
\textit{symplectic core} (of $\mathfrak{m}$).
\begin{lem} Fix notation as above.

(1) The subsets $\mathcal{C}(\mathfrak{m})$ of $\mathcal{Z}$ are
all locally closed, with $\overline{\mathcal{C}(\mathfrak{m})} =
\mathcal{V}(\mathcal{P}(\mathfrak{m}))$ for all $\mathfrak{m} \in
\mathcal{Z}$, if and only if $(A) \Longrightarrow (B)$ of Lemma
\ref{Dix-Moeg} holds for $Z_0$.

(2) Each symplectic core is smooth in its closure.

(2) At least when $k = \C$, the stratification (\ref{strat2}) is a
refinement of (\ref{strat1}).
\end{lem}
\begin{proof} (1) Assume that $\mathcal{C}(\mathfrak{m})$
is locally closed with $\overline{\mathcal{C}(\mathfrak{m})} =
\mathcal{V}(\mathcal{P}(\mathfrak{m}))$, for all $\mathfrak{m} \in
\mathcal{Z}$, and let $\mathfrak{n} \in \mathcal{Z}$. Thus
$$ \mathcal{C}(\mathfrak{n}) = \mathcal{V}(\mathcal{P}(\mathfrak{n}))
\setminus \mathcal{V}(K) $$ for an ideal $K$ of $Z_0$ with
$\mathcal{P}(\mathfrak{n}) \subsetneq K.$ But then $K$ is the
intersection  of all the prime Poisson ideals of $Z_0$ which
strictly contain $\mathcal{P}(\mathfrak{n})$, so that $K$ is
Poisson and (\ref{Dix-Moeg})(B) holds for
$\mathcal{P}(\mathfrak{n})$, as required.

Conversely, assume that (\ref{Dix-Moeg})$(A) \Longrightarrow (B)$
holds for $Z_0$, and let $\mathfrak{n} \in \mathcal{Z}.$ Let $K$
be the intersection of $Z_0$ and the prime Poisson ideals of $Z_0$
which strictly contain $\mathcal{P}(\mathfrak{n})$, so, by
hypothesis, $\mathcal{P}(\mathfrak{n}) \subsetneq K.$ Clearly,
$$ \mathcal{C}(\mathfrak{n}) = \mathcal{V}(\mathcal{P}(\mathfrak{n}))
\setminus \mathcal{V}(K), $$ so that $ \mathcal{C}(\mathfrak{n})$
is locally closed, and the proof is complete.

(2) Let $J$ be the defining ideal of the singular locus of
$\overline{\mathcal{C}(\mathfrak{m})}$. Since $J$ is a Poisson
ideal properly containing $\mathcal{P}(\mathfrak{m})$, it follows
from the definition of a Poisson core that $\mathfrak{m}$ does not
contain $J$. Hence the point corresponding to $\mathfrak{m}$
belongs to the smooth locus of
$\overline{\mathcal{C}(\mathfrak{m})}$, as required.

(3) If $\mathfrak{n} \in \mathcal{C}(\mathfrak{m})$ is a point at
which the rank of $\{-,-\}$ is minimal in
$\mathcal{C}(\mathfrak{m})$, say $\text{rk}(\mathfrak{n}) = r$,
then by Lemma \ref{rank}.2, writing $J_r$ for the defining ideal
of $\mathcal{Z}_r$, $J_r \subseteq \mathcal{P}(\mathfrak{n}) =
\mathcal{P}(\mathfrak{m})$, so that $\text{rk}(\mathfrak{a}) \leq
r$ for all $\mathfrak{a} \in
 \mathcal{C}(\mathfrak{m})$. Then minimality of $r$ implies that the rank is
 constant across  $\mathcal{C}(\mathfrak{m})$.
\end{proof}

The example of the trivial Poisson bracket on $Z_0 = \C [X]$ shows
that in general (\ref{strat2}) is finer than (\ref{strat1}).

\subsection{The Dixmier-Moeglin equivalence revisited}
\label{D-M2} Here is one case in which Question \ref{Dix-Moeg} has
a positive answer.

\begin{lem} Suppose that $Z_0$ has only finitely many Poisson
primitive ideals. Then $\mathcal{P}-\text{spec} (Z_0)$ is a finite
set and (\ref{Dix-Moeg})(A), (B) and (D) are equivalent.
\end{lem}
\begin{proof} With the stated hypothesis, $\mathcal{P}-\text{spec} (Z_0)$ is a finite set by (\ref{P-nullstell}). Let $P \in \mathcal{P}-\text{spec}
(Z_0).$ Then $P = \cap \{ \mathfrak{m} : \mathfrak{m} \in
\mathcal{Z}, P \subseteq \mathfrak{m} \},$ so $P$ is the
intersection of the Poisson cores of the maximal ideals which
contain it. But by hypothesis this intersection is finite, and so
since $P$ is prime we have $P = \mathcal{P}(\mathfrak{m})$ for
some maximal ideal $\mathfrak{m}.$ The equivalence of
(\ref{Dix-Moeg})(A), (B) and (D) now follows from Lemma
\ref{Dix-Moeg}.
\end{proof}

\noindent \textbf{Remark:} In unpublished work, Goodearl obtains a
considerable improvement of the above lemma -- he shows that
Question \ref{Dix-Moeg} has a positive answer provided $Z_0$
admits a torus $H = (k^{\times})^r$ of $k$-algebra Poisson
automorphisms acting rationally, such that
$\mathcal{P}-\text{spec}(Z_0)$ has only finitely many $H$-orbits
of Poisson primitive ideals.

\subsection{The stratification by symplectic leaves}
\label{leaves} For the rest of the paper the underlying base field
$k$ will be $\mathbb{C}$.

We assume initially that $Z_0$ is smooth. By considering the
complex topology, $\mathcal{Z}$ is a complex analytic manifold.
Let $\hat{Z}_0$ be the ring of complex analytic functions on
$\mathcal{Z}$. Then $\hat{Z}_0$ is a Poisson algebra whose bracket
we also denote by $\{ - , - \}$. Indeed, the existence of a
Poisson bracket is equivalent to the existence of a closed
algebraic 2-form on $\mathcal{Z}$; considering the given 2-form as
complex analytic yields the Poisson bracket. Recall that the
\textit{symplectic leaf} $\mathcal{L}(\mathfrak{m})$ containing
the point $\mathfrak{m}$ of $\mathcal{Z}$ is the maximal connected
complex analytic manifold in $\mathcal{Z}$ such that $\mathfrak{m}
\in \mathcal{L}(\mathfrak{m})$ and $\{-,-\}$ is nondegenerate at
every point of $\mathcal{L}(\mathfrak{m})$. By \cite[Proposition
1.3]{wein} symplectic leaves exist and afford a stratification of
$\mathcal{Z}$. By construction the points of the symplectic leaf
$\mathcal{L}(\mathfrak{m})$ are exactly those which can be reached
by travelling along the integral flows of Hamiltonian vector
fields $\{H, -\}$ for $H\in\hat{Z}_0$.

We now extend this definition to an arbitrary commutative Poisson
$\C$-algebra $Z_0$. So, drop the hypothesis that $Z_0$ is smooth,
and define inductively an ascending chain of ideals of $Z_0$, as
follows: set $I_0 = \sqrt \{ 0 \}$ and let $I_{t+1}$ be the ideal
of $Z_0$ such that $I_{t+1}/I_t$ defines the singular locus of
$Z_0/I_t$. Fix $m$ such that $I_m = Z_0.$ By \cite[Proposition
15.2.14(i)]{McC-Rob} or \cite[Corollary 2.4]{Pol} each $I_t$ is a
semiprime Poisson ideal of $Z_0$. Set $Z_t = Z/I_t$ and $A_t =
A/I_tA$. In view of (\ref{g1}) there is an inclusion $Z_t\subseteq
A_t$ on which $D$ induces a map such that the assumptions of
(\ref{setup}) are satisfied.

Let $\mathcal{Z}_t = \text{Maxspec}(Z_t)$ and let
$(\mathcal{Z}_t)_{\text{sm}}$ be the smooth locus of
$\mathcal{Z}_t$. As $(\mathcal{Z}_t)_{\text{sm}}$ is a complex
analytic Poisson manifold it has, as above, a foliation by
symplectic leaves: for each $t = 0, \ldots , m,$ there is an index
set $\mathcal{I}_t$ such that
\[
(\mathcal{Z}_t)_{\text{sm}} = \coprod_{i\in \mathcal{I}_t}
\mathcal{S}_{t,i}.
\]
Since $\mathcal{Z}$ is the disjoint union of the subsets
$(\mathcal{Z}_t)_{\text{sm}}$, there is a stratification of
$\mathcal{Z}$,
\begin{equation}
\label{strat3}  \mathcal{Z} = \coprod_{0 \leq t \leq m,i\in
\mathcal{I}_t}\mathcal{S}_{t,i}.
\end{equation}

We write $\mathcal{L}(\mathfrak{m})$ for the leaf containing
$\mathfrak{m} \in \mathcal{Z}$, and $\overline{\mathcal{S}}_{t,i}$
for the Zariski closure of $\mathcal{S}_{t,i}$. Let $K_{t,i}$ be
the defining ideal of $\overline{\mathcal{S}}_{t,i}$.
\begin{lem} For all $t$, $0 \leq t \leq m,$ and all $i \in
\mathcal{I}_t,$ $K_{t,i}$ is prime and Poisson. Indeed $K_{t,i}$
is the Poisson core $\mathcal{P}(\mathfrak{m})$ of every
$\mathfrak{m} \in \mathcal{S}_{t,i}$.
\end{lem}
\begin{proof}
Let $\mathfrak{m}$ be in $\mathcal{S}_{t,i}$. We show first that
\begin{equation}\label{include1} \mathcal{P}(\mathfrak{m}) \subseteq K_{t,i}.
\end{equation}
So, write $\bar{\mathcal{Z}} =
\text{maxspec}(Z_0/\mathcal{P}(\mathfrak{m}))$, and use
$\bar{\:\:}$ to denote images in
$\hat{Z}_0/\mathcal{P}(\mathfrak{m})\hat{Z}_0$ below. Let
$B(\epsilon)$ be the complex analytic ball with radius $\epsilon$.
Let $H \in \hat{Z}_0,$ and let $\sigma_{\mathfrak{m}}(z) :
B(\epsilon) \longrightarrow \mathcal{Z}$ and
$\overline{\sigma}_{\mathfrak{m}}(z) : B(\epsilon) \longrightarrow
\overline{\mathcal{Z}}$ be the integral curves to $\{H,-\}$ and
$\{\overline{H}, -\}$ respectively, with
$\sigma_{\mathfrak{m}}(0)=\mathfrak{m}$ and
$\overline{\sigma}_{\mathfrak{m}}(0)= \overline{\mathfrak{m}}$.
Thinking of $\overline{\mathcal{Z}}$ as a subset of $\mathcal{Z}$,
it makes sense to claim, as we do, that
\begin{equation} \label{curve} \overline{\sigma}_{\mathfrak{m}} =
\sigma_{\mathfrak{m}}
\end{equation}
in a neighbourhood of $0$. To prove (\ref{curve}), let $f \in
\hat{Z}_0$. By definition of an integral curve,
\begin{equation} \label{der1} \frac{d}{dz} (f \circ \sigma_{\mathfrak{m}}) = \{H,f\} \circ
\sigma_{\mathfrak{m}},
\end{equation}
and
\begin{equation} \label{der2} \frac{d}{dz} (\overline{f} \circ
\overline{\sigma}_{\mathfrak{m}}) = \{\overline{H},\overline{f}\}
\circ \overline{\sigma}_{\mathfrak{m}}.
\end{equation}
But the left hand side of (\ref{der2}) is $\frac{d}{dz} (f \circ
\overline{\sigma}_{\mathfrak{m}})$, while, because
$\mathcal{P}(\mathfrak{m})\hat{Z}_0$ is a Poisson ideal, the right
hand side is $\{H,f\} \circ \overline{\sigma}_{\mathfrak{m}}$.
Hence, comparing this with (\ref{der1}), we deduce from the
uniqueness of flows that $\sigma_{\mathfrak{m}} =
\overline{\sigma}_{\mathfrak{m}}$ in a neighbourhood of $0$. Since
symplectic leaves are  by definition obtained by travelling along
integral curves to Hamiltonians, the leaf through $\mathfrak{m}$
is contained in $\mathcal{V}(\mathcal{P}(\mathfrak{m}))$ and so
(\ref{include1}) follows.

Now we complete the proof by showing that $K_{t,i}$ is a Poisson
ideal. In doing so we may assume without loss that $\mathcal{Z}$
is smooth, since, as in the definition of symplectic leaves above,
we pass to a smooth variety containing $\mathfrak{m}$ when
determining $\mathcal{S}_{t,i}$. So let $f \in K_{t,i}$ and let $H
\in Z_0$. Let $\sigma_{\mathfrak{m}}(z) : B(\epsilon)
\longrightarrow \mathcal{Z}$ be an integral curve to $\{H,-\}$
with $\sigma_{\mathfrak{m}}(0) = \mathfrak{m}$. Then, by
definition of an integral curve, (\ref{der1}) holds. On a complex
analytic neighbourhood of $\mathfrak{m}$, $f \circ
\sigma_{\mathfrak{m}}= 0$ since $\sigma_{\mathfrak{m}}$ has image
in $\mathcal{S}_{t,i}$. Hence
$$ \frac{d}{dz}(f \circ \sigma_{\mathfrak{m}}) = 0 = \{H,f\} \circ
\sigma_{\mathfrak{m}} (0) = \{H,f\}(\mathfrak{m}).$$ That is,
$\{Z_0,K_{t,i}\} \subseteq  \mathfrak{m}$. Repeating this argument
with $\mathfrak{m}$ replaced by each of the members of
$\mathcal{L}(\mathfrak{m})$ in turn we conclude that
$\{Z_0,K_{t,i}\} \subseteq K_{t,i}$ as required.
\end{proof}
\subsection{}\label{alglf}
Of particular interest are the cases where, for all $t$ and for
all $i\in \mathcal{I}_t,$
\begin{equation*}
\label{h1} \mathcal{S}_{t,i} \textit{ is a locally closed
subvariety of } \mathcal{Z}.
\end{equation*}
When this is so we shall say that the Poisson bracket $\{-,-\}$ is
\textit{algebraic}.


The following proposition relates the three stratifications
discussed.
\begin{prop} Let $\mathfrak{m} \in \mathcal{Z}$ with
$\mathrm{rk}(\mathfrak{m}) = j.$
\begin{enumerate}
\item Then
$$ \mathcal{L}(\mathfrak{m}) \subseteq \mathcal{C}(\mathfrak{m}) \subseteq
\mathcal{Z}^o_j.$$ \item Suppose that the Poisson bracket is
algebraic. Then the first inclusion in (1) is an equality, so that
the symplectic leaf of each maximal ideal is then determined by
its core. Moreover (A) and (B) of (\ref{Dix-Moeg}) are equivalent.
$\mathcal{L}(\mathfrak{m})$ is contained in the smooth locus of
$\overline{\mathcal{L}(\mathfrak{m})}$, and
$$ \dim \mathcal{L}(\mathfrak{m}) = \dim
\overline{\mathcal{L}(\mathfrak{m})} = \mathrm{Krull} \: \mathrm{
dim} \big(Z_0/\mathcal{P}(\mathfrak{m})\big) = j. $$
\end{enumerate} \end{prop}

\begin{proof} (1) The first inclusion is a restatement of part of
the above lemma, and the second is Lemma \ref{cores}.2.

(2) Suppose that the Poisson bracket is algebraic. Then, by the
lemma, the defining ideal of the Zariski closure of the leaf
$\mathcal{L}(\mathfrak{m})$ is the core
$\mathcal{P}(\mathfrak{m})$ of $\mathfrak{m}$. Since this
statement applies also to the leaves in the boundary of
$\overline{\mathcal{L}(\mathfrak{m})}$, we deduce that in this
case $\mathcal{L}(\mathfrak{m}) = \mathcal{C}(\mathfrak{m})$. The
equivalence of (\ref{Dix-Moeg})(A) and (B) now follows from Lemmas
\ref{Dix-Moeg} and \ref{cores}(1). That $\dim
\overline{\mathcal{L}(\mathfrak{m})} = \text{Krull dim}
(Z_0/\mathcal{P}(\mathfrak{m}))$ is now clear because
$\mathcal{P}(\mathfrak{m})$ is the defining ideal of
$\overline{\mathcal{L}(\mathfrak{m})}$, by Lemma \ref{leaves}.

Since $\mathcal{L}(\mathfrak{m})$ is open in
$\overline{\mathcal{L}(\mathfrak{m})}$,
$\mathcal{L}(\mathfrak{m})$ contains a point, say $\mathfrak{q}$,
which is smooth in $\overline{\mathcal{L}(\mathfrak{m})}$. By
definition of $\mathcal{L}(\mathfrak{m})$ the Poisson bracket is
nondegenerate on $\mathfrak{q}/\mathfrak{q}^2 +
\mathcal{P}(\mathfrak{m})$, so that
$$ j = \text{rk}(\mathfrak{q}) =
\dim_{\C}\big(\mathfrak{q}/\mathfrak{q}^2 +
\mathcal{P}(\mathfrak{m})\big) = \text{Krull dim}\big(
Z_0/\mathcal{P}(\mathfrak{m})\big), $$ the last equality holding
thanks to smoothness of $\mathfrak{q}$ in
$\overline{\mathcal{L}(\mathfrak{m})}$. Finally, to see that
$\mathcal{L}(\mathfrak{m})$ is contained in the smooth locus of
$\overline{\mathcal{L}(\mathfrak{m})}$, let $\mathfrak{n}$ be any
point of $\mathcal{L}(\mathfrak{m})$. Nondegeneracy of the Poisson
bracket on $\mathcal{L}(\mathfrak{m}) \subseteq \mathcal{Z}^0_j$
forces
$$ \dim_{\C}\big(\mathfrak{n}/\mathfrak{n}^2 +
\mathcal{P}(\mathfrak{m})\big) = j, $$ and since this integer is
$\text{Krull dim} \big( Z_0/ \mathcal{P}(\mathfrak{m})\big)$,
smoothness follows.
\end{proof}

\noindent \textbf{Remarks:} (1) The following example shows that
symplectic leaves need not be algebraic, \cite[Example 2.37]{van}.
 Let
$Z=\mathbb{C}[x,y,z]$ with Poisson bracket given by $\{ x,y \}=0,
\{x,z\} = \alpha x$ and $\{y,z\} = y$ for some $\alpha \in
\mathbb{R}$. A calculation shows that the symplectic leaves fall
into two families: the single points $(0,0,c)$ for $c\in
\mathbb{C}$ and the sets cut out by the equation $xy^{\alpha} =
C$, where $C\in \mathbb{C}$. If $\alpha=p/q\in \mathbb{Q}$ then
the leaves are algebraic, being described by the equation
$x^qy^p=C$. However, if $\alpha \notin \mathbb{Q}$ then the leaves
are not algebraic varieties. In the second case the symplectic
cores are locally closed: they are the points $(0,0,c)$, for $c\in
\mathbb{C}$, and $(\mathbb{C^*})^2 \times \mathbb{C}$.

(2) The inclusion of $\mathcal{L}(\mathfrak{m})$ in the smooth
locus of $\overline{\mathcal{L}(\mathfrak{m})}$ is in general
strict: consider the example in Remark \ref{Dix-Moeg}. Let
$\mathfrak{m} = \langle x-\alpha, y-\beta, z-\gamma\rangle$ with
$\alpha\beta \neq 0$. Then $\mathcal{P}(\mathfrak{m}) = \langle
x-\beta^{-1}\alpha y\rangle$, so that
$\overline{\mathcal{L}(\mathfrak{m})} = \mathcal{V}(\langle
x-\beta^{-1}\alpha y\rangle)$, which is smooth; whereas
$\mathcal{L}(\mathfrak{m}) = \mathcal{V}(\langle
x-\beta^{-1}\alpha y\rangle)\setminus \langle x,y\rangle$.
\subsection{The case of finitely many leaves}\label{fteleaves}In
later applications to symplectic reflection algebras the set of
symplectic leaves will be finite. We now explore the ramifications
of this hypothesis.

\begin{prop} Continue with the notation and hypotheses of
(\ref{leaves}). Suppose that the stratification (\ref{strat3}) of
$\mathcal{Z}$ into symplectic leaves is finite.
\begin{enumerate}
\item The Poisson bracket is algebraic.
\item Let $\mathfrak{m} \in \mathcal{Z}$ with
$\text{rk}(\mathfrak{m}) = j.$ Then the following subsets of
$\mathcal{Z}$ coincide:
\begin{enumerate}
\item the symplectic leaf $\mathcal{L}(\mathfrak{m})$ containing
$\mathfrak{m}$;
\item $\mathcal{C}(\mathfrak{m})$;
\item the irreducible component of $\mathcal{Z}^o_j$ containing
$\mathfrak{m}$;
\item the smooth locus of the irreducible component of
$\mathcal{Z}_j$ containing $\mathfrak{m}.$
\end{enumerate}
\end{enumerate} \end{prop}

\begin{proof} (1) We argue by noetherian induction. So we can
assume that $\mathcal{Z}$ is irreducible by (\ref{pcl}). Since the
singular locus of $\mathcal{Z}$ is a union of leaves by the
construction in (\ref{leaves}), these leaves are locally closed by
our induction assumption, and it remains to consider the leaves
\begin{equation}\label{list} \mathcal{L}_1,\mathcal{L}_2, \ldots , \mathcal{L}_t
\end{equation}
whose disjoint union equals the smooth locus $\mathcal{S}$ of
$\mathcal{Z}$. By Lemma \ref{leaves}, we can apply the induction
hypothesis to any leaf $\mathcal{L}_i$ whose closure is strictly
contained in $\mathcal{Z}$, to see that such leaves are algebraic.
Removing any such leaves from the list (\ref{list}), it follows
that $ \cup_{i=1}^t \mathcal{L}_i$ is a non-empty open set of
smooth points of $\mathcal{Z},$ with $\overline{\mathcal{L}}_i =
\mathcal{Z}$ for all $i$.

Choose $\mathfrak{m} \in  \cup_{i=1}^t \mathcal{L}_i$ for which
$\mathrm{rk}(\mathfrak{m})$ is maximal, say
$\mathrm{rk}(\mathfrak{m})= j$ and $\mathfrak{m} \in
\mathcal{L}_1.$ Since $\overline{\mathcal{L}}_1 = \mathcal{Z}$,
the closed set $\mathcal{Z}_{j-1}$ is proper in $\mathcal{Z}.$ So,
for all $i = 1, \ldots , t,$ $\overline{\mathcal{L}}_i =
\mathcal{Z} \nsubseteq \mathcal{Z}_{j-1}$, implying that
$\mathcal{L}_i \nsubseteq \mathcal{Z}_{j-1}.$ It follows from this
and the maximality of $\mathrm{rk}(\mathfrak{m})$ that the rank is
constant and equal to $j$ across $ \cup_{i=1}^t \mathcal{L}_i$.
Therefore $ \cup_{i=1}^t \mathcal{L}_i$ is a $j$-dimensional
irreducible smooth open subset of $\mathcal{Z}$ on which the rank
is constant and equal to $j$. Hence, $ \cup_{i=1}^t \mathcal{L}_i$
is a leaf. That is, $t = 1$ and $\mathcal{L}_1$ is locally closed,
as required.

(2) First, (a) = (b) and (b) is contained in (c) by the first part
of the proposition and Proposition \ref{alglf}. Now let
$I(\mathcal{A})$ be the defining ideal of the irreducible
component $\mathcal{A}$ of $\mathcal{Z}^o_j$ containing
$\mathfrak{m}$. We claim that
\begin{equation}\label{kdim} \text{Krull dim}\big(Z_0/I(\mathcal{A})\big) = j.
\end{equation}
For certainly $ \text{Krull dim}\big(Z_0/I(\mathcal{A})\big) \geq
j $ by Lemma \ref{rank}.5. Suppose that the inequality is strict.
Then $\mathcal{A}$ is a union of symplectic leaves by Proposition
\ref{alglf}.1, and these leaves all have closures of dimension $j$
by the first part of the proposition and Proposition
\ref{alglf}.2. Hence infinitely many such leaves are needed to
cover $\mathcal{A}$, contradicting the finiteness hypothesis. So
(\ref{kdim}) is proved.

Now $\mathcal{A}$ is a union of $j$-dimensional algebraic leaves,
but (\ref{kdim}) and the irreducibility of $\mathcal{A}$ means
that the union consists of a single leaf, proving (b) = (c).
Finally, Lemma \ref{rank}.6 shows that (c) is contained in (d). To
prove equality we adapt an idea from \cite[Lemma 3.5(ii)]{EG}. Let
$\mathcal{S}$ be the smooth locus of $\overline{\mathcal{A}}$, and
suppose for a contradiction that $\mathcal{A} \subsetneq
\mathcal{S}$. Now $\mathcal{S}$ is an (algebraic) Poisson
manifold, for which the subvariety $\mathcal{U}$ of points of rank
less than $j$ is non-empty and has codimension at least 2 in
$\mathcal{S}$, since there are by hypothesis only a finite number
of leaves and each one has even dimension. On the other hand
$\mathcal{U}$ has codimension 1 in $\mathcal{S}$ since it is
defined by the vanishing of the divisor associated with the
Poisson form on $\mathcal{S}$. This contradiction shows that (c) =
(d).
\end{proof}

\subsection{G-equivariant stratifications} \label{Gstrat}
Let $G$ be an algebraic group. Suppose that $G$ acts rationally by
algebra automorphisms on $A$, that this action preserves $Z_0$ and
that $D$ is $G$-equivariant, in the sense that $D_{gz}(a) =
gD_z(g^{-1}a)$ for all $g \in G, z \in Z_0$ and $a \in A$. Then
$G$ consists of Poisson automorphisms of $Z_0$.

One can carry through a $G$-equivariant version of most of
paragraphs (\ref{Dix-Moeg}) -- (\ref{fteleaves}), defining, for
example, prime $G$-Poisson ideals, $G$-Poisson primitive ideals,
$G$-smyplectic cores, $G$-symplectic leaves (meaning $G$-orbits of
symplectic leaves) and so on, and proving $G$-equivariant versions
of most of the results of those paragraphs. In discussing the
algebraicity of $G$--symplectic leaves, without loss of generality
we may assume that $G$ is connected, if necessary replacing $G$ by
its identity component $G^o$. We label the $G$--orbits of
symplectic leaves $\mathcal{G}_i$, so that we have a
stratification
\begin{equation} \label{gstr} \mathcal{Z} = \coprod \mathcal{G}_i.
\end{equation}
\begin{prop}
Assume that $G$ is connected.  If the Poisson bracket on $Z_0$ is
algebraic, then each stratum $\mathcal{G}_i$ is locally closed and
irreducible.
\end{prop}
\begin{proof}
Consider the action morphism $\alpha: G \times \mathcal{Z}
\longrightarrow \mathcal{Z}.$ The variety $\mathcal{G}_i$ is
irreducible since it is the image under $\alpha$ of the
irreducible variety $G\times \mathcal{L}$ for some leaf
$\mathcal{L}$. Moreover, $\mathcal{L}$ is locally closed by
hypothesis, and with it so is $G \times \mathcal{L}$. By
Chevalley's theorem, \cite[Exercise II.3.19]{Hartshorne}, the
image under a morphism of a locally closed subvariety is
constructible, so that $\mathcal{G}_i$ is constructible, meaning
(thanks to \cite[Exercise II.3.18]{Hartshorne}) that there is a
finite union
\[
\mathcal{G}_i = \bigcup_s \mathcal{U}_s\cap \mathcal{C}_s,
\]
where $\mathcal{U}_s$ is open in $\mathcal{Z}$ and $\mathcal{C}_s$
is closed in $\mathcal{Z}$ for each $s$. Without loss of
generality we may assume that $\mathcal{U}_s\cap \mathcal{C}_s
\neq \emptyset$ and that $\mathcal{C}_s\subseteq
\overline{\mathcal{G}}_i$ for all $s$. As $\mathcal{G}_i$ is
irreducible, there exists $s'$ such that
$\mathcal{C}_{s'}=\overline{\mathcal{G}}_i$. Hence
$\mathcal{U}_{s'}\cap \mathcal{C}_{s'}$ is open and dense in
$\overline{\mathcal{G}}_i$ and consequently $\mathcal{G}_i$ itself
is dense in $\overline{\mathcal{G}}_i$.

Thus the smooth locus of $\overline{\mathcal{G}}_i$ meets
$\mathcal{G}_i$ and, by applying the $G$-action if necessary, we
see that this smooth locus actually meets the symplectic leaf
$\mathcal{L}$. Since the defining ideals of
$\overline{\mathcal{L}}$ and hence of $\overline{\mathcal{G}}_i$
are Poisson closed by Lemma \ref{leaves}, we can integrate the
derivations arising from $D$ to one-parameter families of complex
analytic automorphisms of
$(\overline{\mathcal{G}}_i)_{\text{sm}}$, acting transitively on
$\mathcal{L}$. Combining this again with the action of $G$ we
deduce
 that each point of
$\mathcal{G}_i$ is smooth in $\overline{\mathcal{G}}_i$. Finally,
$\mathcal{U}_{s'}\cap \mathcal{C}_{s'}$ is open in
$\overline{\mathcal{G}_i}$, so that at least one point of
$\mathcal{G}_i$ has an open neighbourhood
 in $\overline{\mathcal{G}}_i$ which is
contained in $\mathcal{G}_i$. Applying $G$ and the automorphisms
introduced above it follows that every point of $\mathcal{G}_i$
has such a neighbourhood, so that $\mathcal{G}_i$ is open in
$\overline{\mathcal{G}_i}$.
\end{proof}

\section{The bundle of finite dimensional algebras}\label{iss}
In this section we continue to assume that $A$ is a Poisson
$Z_0$-order, and we'll freely use the notation introduced in the
earlier sections. The underlying field $k$ will always be $\C$.
\subsection{} \label{pollem} Let $M$ be a $Z_0$-module. We say that $M$ is a
\textit{Poisson module} if there exists a bilinear form $\{ - , -
\}_M : Z_0\times M \longrightarrow M$ satisfying $\{z,z'm\}_M =
\{z,z'\}m + z'\{z,m\}_M$ for all $z,z'\in Z_0$ and $m\in M$. The
following is proved in \cite[Lemma 2.1]{Pol} \begin{lem} Let $M$
be a finitely generated $Z_0$-Poisson module. Then the annihilator
of $M$ is a Poisson ideal of $Z_0$, and for any $n$ the ideal
defining the locus where the rank of $M$ as a $Z_0$-module is
greater than $n$ (in the sense of the minimal number of
generators) is Poisson closed.
\end{lem}
\begin{proof}
Let $z'M=0$ for some $z'\in Z_0$. Then for all $z\in Z_0$ and
$m\in M$
\[
0 = \{z, z'm\}_M = \{z,z'\}m + z'\{z,m\}_M = \{z,z'\}m,
\]
proving that the annihilator of $M$ is a Poisson ideal. The second
claim follows by applying the first to the exterior powers of $M$.
\end{proof}
\subsection{} \label{algIM}
The results we shall prove in this Section were proved by De
Concini, Lyubashenko and Procesi \cite[Corollary 11.8]{DeCP},
\cite[Corollary 9.2]{DeCL}, in the special case of a Poisson
$Z_0$-order $A$ with base field $\C$, such that $Z_0$ has finite
global dimension and $A$ is a free $Z_0$-module, and an algebraic
Poisson bracket. In the Hopf algebra setting of quantum groups
these additional hypotheses are valid, but they no longer hold for
symplectic reflection algebras and we are thus obliged to prove
the stronger results given here.\footnote{We are grateful to
several participants in the Oberwolfach meeting in Noncommutative
Geometry of 14-19 April 2002, and in particular Alastair King, for
helpful comments in conection with Theorem (\ref{algIM}).} Much of
our proof is an adaptation of the earlier ones, but since the
latter were somewhat brief we have included full details here as
an aid to the reader.

\begin{thm}
Let $A$ be a Poisson $Z_0$-order which is an affine $\C$-algebra.
For each point $x\in\mathcal{Z}$, define the finite dimensional
algebra
\[ A_{x} := \frac{A}{\mathfrak{m}_{x}A}.
\]If $x,y$ belong
to the symplectic core $\mathcal{C}$, then $A_{x}$ and $A_{y}$ are
isomorphic $\C$-algebras.
\end{thm}
\begin{proof}
Let $I$ be the defining ideal of $\overline{\mathcal{C}}$, a prime
Poisson ideal of $Z_0$ by Lemma \ref{cores}. Both $A_x$ and $A_y$
are quotients of $A/IA$, an algebra satisfying the hypotheses of
(\ref{setup}). In particular $A/IA$ is a Poisson $Z_0/I$-module
with bracket induced from $D_z$ for $z\in Z_0$. By Lemma
\ref{pollem}, $\mathcal{D}$, the subvariety of points of
$\overline{\mathcal{C}}$ for which the rank of $A/IA$ as a module
over $Z_0/I$ is non-minimal, is Poisson closed and, by definition,
proper in $\overline{\mathcal{C}}$. Set $\mathcal{E} =
(\overline{\mathcal{C}}\setminus \mathcal{D})_{sm}$, a locally
closed variety with closure $\overline{\mathcal{C}}$. By
construction the rank of $A/IA$ as a module over $Z_0/I$ is
constant on $\mathcal{E}$. It follows from \cite[Exercise
II.5.8(c)]{Hartshorne} that $A/IA$ is free (of finite rank) over
$\mathcal{E}$.

The sets $\mathcal{C}$ and $\mathcal{D}$ are disjoint. Indeed, if
$x\in \mathcal{C}\cap \mathcal{D}$ then, since $\mathcal{C}$ is
the Poisson core of $x$, $\overline{\mathcal{C}} \subseteq
\mathcal{D}$, contradicting properness. It follows from Lemma
\ref{cores}(2) that $\mathcal{C} \subseteq \mathcal{E}$.

We will now argue as in \cite[Section 9]{DeCL}. Thanks to the
above,  we may assume without loss of generality that $\mathcal{E}
= \mathcal{Z}$ is smooth and $A$ is a free $Z_0$-module. Recall
that $\hat{Z_0}$ denotes the ring of complex analytic functions on
$\mathcal{Z}$. Define $\hat{A} = A\otimes_{Z_0} \hat{Z}_0$. Note
that, for any value of $x$, the natural map $$A_x =
\frac{A}{\mathfrak{m}_xA} \longrightarrow
\frac{\hat{A}}{\hat{\mathfrak{m}}_x\hat{A}} = \hat{A}_x$$ is an
algebra isomorphism. Moreover, given any $H\in Z_0$, then the
derivation $D_H$ extends uniquely to a derivation, $\hat{D}_H$, on
$\hat{A}$, thanks to the extension of the Poisson bracket on $Z_0$
to $\hat{Z}_0$.

Consider $\hat{A}$ as a vector bundle $\mathcal{V}$ over
$\mathcal{Z}$ with fibres $\hat{A}_x$. Since this bundle is
trivial, we have $\mathcal{V}\cong \mathcal{Z} \times V$ where $V$
is a vector space isomorphic to $A_x$. Pick a basis for $A$ over
$Z_0$, say $\{a_1,a_2, \ldots ,a_n\}$ with $a_1=1$. The above
isomorphism identifies $\sum\lambda_ia_i \in \hat{A}_x$ with
$(x,(\lambda_1,\ldots ,\lambda_n))$.

To the derivation $\hat{D}_H$ we can associate a vector field,
$\xi_H$, on $\mathcal{V}$. Since $T\mathcal{V} =
T(\mathcal{Z}\times V)$, if $y= (x, (\lambda_1,\ldots
,\lambda_n))$ we have $T_y\mathcal{V} = T_x\mathcal{Z} \times
\mathbb{C}^n$, and we write the vector field as
\[
{\xi}_{H,y} = \left(\xi_{H,x}, \left(-\sum_{i=1}^n
p_H^{i,j}(x)\right)_j\right),
\]
where $\xi_{H,x}(f) = \{ H,f\} (x)$ is the Hamiltonian vector
field on $\mathcal{Z}$ evaluated at $x$ and $\{ H,a_i\}=
\sum_{j=1}^n p_H^{i,j}a_j$ with $p_H^{i,j}\in \hat{Z}_0$.

Let $\phi : B(\epsilon)\times \mathcal{V} \longrightarrow
\mathcal{V}$ be a local flow around $y$ lifting $\xi_H$. Note that
$\phi_z$ sends fibres to fibres since the Hamiltonian vector field
of $\xi_H$ is independent of $\lambda_1,\ldots ,\lambda_n$. It
follows from the definition of $\xi_{H,y}$ that $\phi_z$ is linear
on the fibres also.

Let $\varrho: B(\epsilon) \times \mathcal{Z} \longrightarrow
\mathcal{Z}$ be the local flow on $\mathcal{Z}$ around $x$,
induced by the Hamiltonian $\xi_{H}$. Splitting $\phi$ into
components, we have a linear isomorphism $$\psi_z:
\hat{A}_{\varrho(0)} \longrightarrow \hat{A}_{\varrho(z)},$$ which
is given explicitly in \cite[9.1]{DeCL}. We claim that $\psi_z$ is
an isomorphism of algebras. Up to first order, $\psi_z$ is
described as $$\text{id}-\hat{D}_H : \hat{A}_{\varrho(z)}
\longrightarrow \hat{A}_{\varrho(z+dz)}.$$ To prove this is an
algebra isomorphism, let multiplication in $\hat{A}_x$ be denoted
by $\mu_x$ and $c_{ij}^k\in \hat{Z}_0$ be the structure constants
of $\hat{A}$ with respect to the chosen basis $\{a_1, \ldots , a_n
\}$ . We have
\begin{eqnarray*}
&&\mu_{\varrho(z+dz)}((\text{id}-\hat{D}_{H,\varrho(z)}dz)(a_i)\otimes
(\text{id}-\hat{D}_{H,\varrho(z)}dz)(a_j))= \\ && \qquad \qquad
=\mu_{\varrho(z+dz)} ((a_i -\{H,a_i\}dz )\otimes (a_j -
\{H,a_j\}dz)) \\ && \qquad \qquad=\sum_{k=1}^n c_{ij}^k
(\varrho(z+dz))a_k - \mu_{\varrho(z)}(\{H,a_i\}\otimes a_j +
 a_i\otimes \{H,a_j\})dz \\
&& \qquad \qquad=\left(\sum_{k=1}^n c_{ij}^k(\varrho(z))a_k +
\sum_{k=1}^n \{H,c_{ij}^k\}(\varrho(z))dz \right) -
\{H,a_ia_j\}dz \\
&& \qquad \qquad= \sum_{k=1}^n c_{ij}^k(\varrho(z))a_k + \left(
\sum_{k=1}^n \{H,c_{ij}^k\}(\varrho(z)) -\sum_{k=1}^n(
\{H,c_{ij}^k\}(\varrho(z)) + c_{ij}^k(\varrho(z))\{H,
a_k\})\right)dz
\\&&  \qquad \qquad=(\text{id}-\hat{D}_H)\mu_{\varrho(z)}(a_i\otimes a_j).
\end{eqnarray*}
Let $d_{ij}^k\in \hat{Z}_0$ be the structure constants obtained by
transport of structure along $\psi_z$. Both $c_{ij}^k$ and
$d_{ij}^k$ provide analytic maps $B(\epsilon)\longrightarrow
\text{Alg}_{\mathbb{C}}(n)$, where $\text{Alg}_{\mathbb{C}}(n)$
denotes the variety of $n$-dimensional algebras over $\mathbb{C}$.
The above calculation shows that the derivatives of these maps
agree everywhere. Therefore the maps are equal, proving the claim
that $\psi_z$ is an algebra isomorphism.

We have now shown that any Hamiltonian flow generated by an
algebraic function lifts to an isomorphism of
$\mathbb{C}$-algebras within the leaf $\mathcal{L}_x$. Since the
set of algebraic functions is dense in the set of analytic
functions it follows that we can trace out a dense subset of
$\mathcal{L}_x$ by algebraically generated Hamiltonian flows. Call
this set $\mathcal{L}^{\text{alg}}_x$. By Lemma \ref{leaves} we
have
$$\overline{\mathcal{L}^{\text{alg}}_x}\cap \mathcal{E} =
\overline{\mathcal{L}_x}\cap \mathcal{E} = \mathcal{E}.$$

There is an action of $GL_n(\mathbb{C})$ on
$\text{Alg}_{\mathbb{C}}(n)$ by base change, whose orbits are the
isomorphism classes of $n$--dimensional algebras over
$\mathbb{C}$. Let $\Upsilon : \mathcal{E} \longrightarrow
\text{Alg}_{\mathbb{C}}(n)$ be the morphism obtained by sending a
point $x\in \mathcal{E}$ to $A_x$. By the above paragraph
$\Upsilon(\mathcal{L}^{\text{alg}}_x)$ is contained in a unique
$GL_n(\mathbb{C})$--orbit, say $\mathcal{O}_x$. As
$\mathcal{L}^{\text{alg}}_x$ is dense in $\mathcal{E}$, the image
of $\Upsilon$ lies in $\overline{\mathcal{O}_x}$. Repeating the
above argument for $y\in \mathcal{C}\subseteq \mathcal{E}$, we
deduce that $\overline{\mathcal{O}_x} = \overline{\mathcal{O}_y}$.
Therefore $\mathcal{O}_x = \mathcal{O}_y$, showing that $A_x\cong
A_y$ as required.
\end{proof}

\noindent \textbf{Remark:} In situations where Question
\ref{Dix-Moeg} has a positive answer, Lemma \ref{cores}(1) shows
that the symplectic cores are locally closed in $\mathcal{Z}$. In
this situation, the symplectic reflection algebra example in
Section \ref{difffd} (taking $\Gamma = \mathbb{Z}_2$) shows that
the theorem cannot be improved: different symplectic cores can
yield non--isomorphic algebras.
\subsection{$G$-equivariant isomorphisms}\label{minstrat}
Let $G$ be a connected algebraic group, acting on $A$ as in
(\ref{Gstrat}). We have a stratification by $G$--symplectic cores,
\[ \mathcal{Z} = \coprod \mathcal{GC}_i.
\] Clearly
$\mathfrak{p}_i = I(\overline{\mathcal{GC}_i})$ is a Poisson
closed prime ideal
 of $Z_0$. Hence, by  (\ref{pcl}), both $\mathfrak{p}_iA$ and the minimal
primes lying over it are Poisson ideals of $A$.

\begin{prop} Retain the hypotheses of (\ref{Gstrat}) and the notation just introduced. Let $P$ be a prime
ideal of $A$ minimal over $\mathfrak{p}_iA$. For all
$\mathfrak{m}_x,\mathfrak{m}_y\in\mathcal{GC}_i$ there are algebra
isomorphisms
\[
A/\mathfrak{m}_xA \cong A/\mathfrak{m}_yA \quad \text{ and } \quad
\frac{A}{P+\mathfrak{m}_xA} \cong \frac{A}{P+\mathfrak{m}_yA}.
\]
\end{prop}
\begin{proof}
Since $G$ is connected and preserves $\mathfrak{p}A$, $g (P) = P$
for all $g \in G.$ Thus, by suitable application of elements of
$G$, we may assume that $\mathfrak{m}_x$ and $\mathfrak{m}_y$ are
in the same symplectic core. It is now clear that the algebra
isomorphism from $A_x$ to $A_y$ afforded by Proposition
\ref{algIM} maps $P$ to itself. The result follows.
\end{proof}

\section{Azumaya sheaves}
\label{az}
\subsection{}\label{setup2}
In this section we examine in more detail the conclusion derived
in Proposition \ref{minstrat}. Thus we impose the following
hypotheses throughout Section \ref{az}:
\begin{enumerate}
\item $A$ is a noetherian $\C$-algebra, finitely generated over a central
subalgebra $Z_0$;
\item Let $\mathcal{Z} = \text{Maxspec}Z_0$. There is a stratification
\[
\mathcal{Z} = \coprod \mathcal{Z}_i,
\]
where each $\mathcal{Z}_i$ is an irreducible locally closed
subvariety of $\mathcal{Z}$ such that $\overline{\mathcal{Z}}_i$,
the closure of $\mathcal{Z}_i$, is a union of some of the $\mathcal{Z}_j$;
\item Let $\mathfrak{p}_i = I(\overline{\mathcal{Z}_i})$ be the prime ideal of $Z_0$ defining
$\overline{\mathcal{Z}_i}$, and let $P_{i,j}$ be the minimal
primes of $A$ lying over $\mathfrak{p}_iA$. For all $x,y \in
\mathcal{Z}_i$ and all $j$ there is an algebra isomorphism
\[
\frac{A}{P_{i,j}+ \mathfrak{m}_xA}\cong
\frac{A}{P_{i,j}+\mathfrak{m}_yA}.
\]
\end{enumerate}
\subsection{}\label{sheaf}
Recall that $\mathcal{A}$ is a sheaf of algebras over a
 variety $\mathcal{V}$ if $\mathcal{A}$ is a quasi-coherent sheaf
 over $\mathcal{V}$ such that for each open set
$\mathcal{U}\subseteq \mathcal{V}$, the sections
$\mathcal{A}(\mathcal{U})$ yield an
$\mathcal{O}_{\mathcal{V}}(\mathcal{U})$-algebra. Moreover, if
there exists an open affine covering of $\mathcal{V}$, say $\{
\mathcal{V}_i\}$, such that $\mathcal{A}(\mathcal{V}_i)$ is an
Azumaya algebra finitely generated as a module over
$\mathcal{O}_{\mathcal{V}}(\mathcal{V}_i)$ for each $i$, then we
shall call $\mathcal{A}$ a \textit{sheaf of Azumaya algebras} over
$\mathcal{V}$. The proof of the following routine lemma is left to
the reader.

\begin{lem} Suppose $\mathcal{A}$ is a sheaf of Azumaya algebras
over a noetherian variety $\mathcal{V}$. Then, for any open
irreducible affine subvariety $\mathcal{U}$, the algebra
$\mathcal{A}(\mathcal{U})$ is Azumaya.
\end{lem}

We come now to the main result of this section.
\begin{prop}
Retain the assumptions of  (\ref{setup2}). Let
$\mathcal{A}_{i,j}$ be the sheaf of algebras over
$\overline{\mathcal{Z}}_{i}$ corresponding to $A/P_{i,j}$. Then
the restriction of $\mathcal{A}_{i,j}$ to $\mathcal{Z}_{i}$ is a
sheaf of Azumaya algebras.
\end{prop}

\begin{proof} By (\ref{setup2})(3), $P_{i,j} \cap Z_0 = \mathfrak{p}_i$, so
we can consider the triple $$Z_0/\mathfrak{p}_i \subseteq
Z(A/P_{i,j}) \subseteq A/P_{i,j}.$$ By (\ref{setup2})(1) the
extension $Z_0/\mathfrak{p}_{i}\subseteq Z(A/P_{i,j})$ is
generically separable by \cite[Lemma 2.1]{BrGor} and the prime
ring $A/P_{i,j}$ is generically Azumaya by \cite[13.7.4 and
13.7.14]{McC-Rob}. Hence there exists a non-empty open set
$\mathcal{U}\subseteq \overline{\mathcal{Z}}_i$ and integers $d$
and $s$, such that for all $x\in\mathcal{U}$ we have
\begin{equation} \label{azum} \frac{A}{P_{i,j}+\mathfrak{m}_xA}
\cong \bigoplus_{i=1}^s \text{Mat}_d(\C). \end{equation} By
(\ref{setup2})(2) $\mathcal{Z}_i$ and $\mathcal{U}$ have non-empty
intersection, and so $\mathcal{Z}_i \subseteq \mathcal{U}$ by
(\ref{setup2})(3).

Let $f\in Z_0/\mathfrak{p}_i$ be any non-zero function vanishing
on $\overline{\mathcal{Z}}_i\setminus \mathcal{Z}_i$. Then
$(\overline{\mathcal{Z}}_i)_{f}\subseteq \mathcal{Z}_i$ and so, by
(\ref{azum}) and the Artin-Procesi theorem,
\cite[13.7.14]{McC-Rob}, $A/P_{i,j}[f^{-1}]$ is an Azumaya
algebra. The proposition follows by covering $\mathcal{Z}_i$ by
such distinguished open sets.
\end{proof}

\section{Quantised function algebras}\label{qfa} \subsection{}Let $G$ be a
simply-connected, semisimple algebraic group over $\C$ and let $T$
be a maximal torus contained in a Borel subgroup $B$ of $G$. Let
$B^{-}$ be the Borel subgroup of $G$ opposite $B$ and let $W$ be
the Weyl group of $G$ with respect to $T$.

Let $A=\EO$ be the quantised function algebra of $G$ at a
primitive $\ell$th root of unity $\epsilon$, where $\ell$ is odd,
and prime to 3 if $G$ contains a factor of type $G_2$. There is a
central subalgebra $Z_0$ of $\EO$ isomorphic to $\CO$, with $\EO$
a finitely generated $Z_0$-module. Using the construction in
(\ref{alg1}), it is shown in \cite{DeCL} that the pair $A$ and
$Z_0$ satisfy the conditions of (\ref{setup}). Moreover, the
symplectic leaves are algebraic, \cite[Appendix A]{H-L}.

\subsection{}There is an action of $T$ on $A$ by winding
automorphisms, preserving $Z_0$ and acting as Poisson
automorphisms of the latter,
 \cite[Propositions 9.3 and 8.7(b)]{DeCL}. Under this action the orbits of the
symplectic leaves are identified in \cite[Appendix A]{H-L} with the double Bruhat cells of
$G$, \[ G = \coprod_{w_1,w_2\in W} X_{w_1,w_2},\] where
$X_{w_1,w_2} = B\dot{w}_1B\cap B^{-}\dot{w}_2B^{-}$, with $\dot{w}_i$
denoting an inverse image of $w_i \in W$ in $N_G(T).$

\subsection{}For arbitrary $g\in G$ the algebras $A_g$ of
(\ref{algIM}) are rather incompletely understood -- for the
current state of knowledge, see \cite{BrGor2}. On the other hand,
for any $w_1,w_2\in W$ the sheaf of Azumaya algebras lying over
$X_{w_1,w_2}$ of Proposition \ref{sheaf} is in this case
explicitly described in \cite{DeCP2}.

\section{Symplectic reflection algebras}\label{sra}
\subsection{Definition and fundamental properties}\label{sradefs} Given the
data of a $2n-$dimensional complex symplectic vector space $V$ and
a finite subgroup $\Gamma$ of $\text{Sp}(V)$, Etingof and Ginzburg
\cite{EG} construct a family of so-called \textit{symplectic
reflection algebras}, as follows. First, define an element $s \in
\Gamma$ to be a {\em symplectic reflection} if (in its action on
$V$) $\text{rank}(1 - s) = 2.$ The set $S$ of symplectic
reflections in $\Gamma$ is closed under conjugation. Take a
complex number $t$ and a $\Gamma$-invariant function ${\bf c}:S
\longrightarrow \C: s \mapsto c_s.$ Note that for $s \in \Gamma$
there is a $\omega$-orthogonal decomposition $V = \mathrm{Im}(1 -
s) \oplus \mathrm{Ker}(1 - s).$ For $s \in S$, write $\omega_s$
for the skew-symmetric form on $V$ which has $\mathrm{Ker}(1 - s)$
as its radical, and coincides with $\omega$ on $\mathrm{Im}(1 -
s)$. Now define $A_{t, \bf{c}}$ to be the $\C$-algebra with
generators $V$ and $\Gamma$, and relations those for $\Gamma$,
together with $$ \gamma x \gamma^{-1} = \gamma (x), \text{ and }
xy - yx = t \omega (x,y)1_{\Gamma} + \sum_{s \in S} c_s \omega_s
(x,y)s, $$ for $x,y \in V$ and $\gamma \in \Gamma$. Thus, for a
given triple $(V,\omega,\Gamma)$ these algebras form a family
parametrised by a complex number $t$ together with the points of
an affine space of dimension $r$ (where $r$ is the number of
conjugacy classes of symplectic reflections in $\Gamma$). It's
clear from the relations above that
\begin{equation} \label{projcong} A_{t,\bf{c}} \cong A_{\mu t, \mu \bf{c}}\end{equation}
 for $\mu \in
\C^*$, so that there is a space $\mathbb{P}^{r}(\C)$ of symplectic
reflection algebras arising from a given $(V,\omega,\Gamma)$. This
space includes the familiar special cases $A_{0, \bf{0}}$ and
$A_{1, \bf{0}}$, the skew group algebras of $\Gamma$ over the
algebra $\mathcal{O}(V)$ and the Weyl algebra $A_n(\C)$.

Clearly, $A_{t, \bf{c}}$ is a filtered $\C$-algebra: namely, set
\begin{equation}
\label{filt}  F_0 = \C \Gamma; \quad F_1 = \C \Gamma +  \C \Gamma
V; \quad \mathrm{and} \quad F_i = (F_1)^i, \quad \mathrm{for}
\quad i \geq 1.
\end{equation} We can form the associated graded ring
$\mathrm{gr}(A_{t, \bf{c}})$ of $A_{t, \bf{c}}$. There is an
obvious epimorphism of algebras $$ \rho : \mathcal{O}(V) \ast
\Gamma \twoheadrightarrow \mathrm{gr}(H_{\kappa}). $$  Etingof and
Ginzburg \cite[Theorem 1.3]{EG} prove the beautiful result that
\begin{equation}
\label{PBW}
 \rho \textit{ is an isomorphism}. \end{equation} This is {\em the
PBW theorem for symplectic reflection algebras}. Notice that we
can immediately conclude from (\ref{PBW}) that $A_{t, \bf{c}}$ is
a deformation of $\mathcal{O}(V) \ast \Gamma$, the family
$A_{t,\bf{c}}$ over $\mathbb{P}^r(\mathbb{C})$ is flat, and (using
standard filtered-graded techniques) deduce that $A_{t, \bf{c}}$
is a prime noetherian $\C$-algebra with good homological
properties.

\subsection{The centre of $A_{t,\bf{c}}$}\label{centre}  Let $e=
|\Gamma|^{-1}\sum_{\gamma\in \Gamma}{\gamma}\in \mathbb{C}\Gamma
\subseteq A_{t,\bf{c}}$, so $e$ is the familiar averaging
idempotent. The filtration (\ref{filt}) induces by intersection a
filtration of $eA_{t,\bf{c}}e$, whose associated graded ring is
clearly $e\mathcal{O}(V) \ast \Gamma e$. It's easy to see that the
latter algebra is isomorphic to the invariant ring
$\mathcal{O}(V)^{\Gamma}$, which is the centre of $\mathcal{O}(V)
\ast \Gamma$. In particular, $eA_{t,\bf{c}}e$ is a (not
necessarily commutative) noetherian integral domain. Let
$Z_{t,\bf{c}}$ be the centre of $A_{t,\bf{c}}$. One proves by a
straightforward adaptation of the proof of \cite[Theorem 3.1]{EG}
that $Z_{t,\bf{c}}$ is isomorphic to the centre of the integral
domain $eA_{t,\bf{c}}e$ via the map $z \mapsto ez$.
\begin{prop}
Retain the above notation.
\begin{enumerate}
\item If $t=0$ then $Z_{0,\bf{c}} \cong eA_{0,\bf{c}}e$ and
$A_{0,\bf{c}}$ is a finitely generated $Z_{0,\bf{c}}$-module for
all $\bf{c}\in\mathbb{A}^r$. \item If $t\neq 0$ then $Z_{t,\bf{c}}
= \mathbb{C}$ for all $\bf{c}\in \mathbb{A}^r$. \end{enumerate}
\end{prop}
\begin{proof} (1) This is \cite[Theorem 3.1]{EG}.

(2) Assume that $t \neq 0$. We work in the algebra
$eA_{t,\bf{c}}e$. As explained above, this algebra has a
filtration $F^0\subseteq F^1\subseteq F^2 \subseteq \cdots$, whose
associated graded ring equals $eA_{0,\bf{0}}e \cong
\mathcal{O}(V)^{\Gamma}$. By \cite[Claim 2.25]{EG},
$eA_{t,\bf{c}}e$ is noncommutative: more precisely, given $u\in
F^m$ and $v\in F^n$, the commutator $[u,v]$ is contained in
$F^{m+n-2}$, and moreover there exist $m,n$, $u$ and $v$ such that
$[u,v]\notin F^{m+n-3}$. Thus, following (\ref{p1}) in
(\ref{alg2}) with $m+n-2$ replacing $m+n-1$, we have a non-trivial
Poisson bracket of degree $-2$ on $\mathcal{O}(V)^{\Gamma}$, say
$B(-, -)$. Since $\Gamma \subseteq \text{Sp}(V)$, the Poisson
bracket on $\mathcal{O}(V)$ defined by the symplectic form on $V$
restricts to $\mathcal{O}(V)^{\Gamma}$ -- that is,
$\{\mathcal{O}(V)^{\Gamma}, \mathcal{O}(V)^{\Gamma} \} \subseteq
\mathcal{O}(V)^{\Gamma}$. It's a consequence of Hartog's theorem
that, up to scalar multiplication, this is the unique bracket of
degree $-2$, \cite[Lemma 2.23]{EG}. We thus have $$ B(-,-) =
\lambda \{ - ,-\}$$ for some $\lambda\neq 0$. The principal symbol
of any element $z$ in the centre of $eA_{t,\bf{c}}e$ belongs to
$\text{Cas}(\mathcal{O}(V)^{\Gamma})$, the algebra of Casimirs on
$\mathcal{O}(V)^{\Gamma}$, since \[ \{\sigma_m(z),\sigma_n(u)\} =
\lambda^{-1}B(\sigma_m(z),\sigma_n(u)) = \lambda^{-1}( [z,u] +
F^{m+n-3}) = 0 + F^{m+n-3}.\] Thus it suffices to prove that
$\text{Cas}(\mathcal{O}(V)^{\Gamma}) =\mathbb{C}$, since the
scalars are the only elements of $eA_{t,\bf{c}}e$ for which the
principal symbol lies in $\mathbb{C}$.

Since $V$ is a symplectic vector space,
$\text{Cas}(\mathcal{O}(V))=\mathbb{C}$ . Let $p\in
\text{Cas}(\mathcal{O}(V)^{\Gamma})$ and let $u\in
\mathcal{O}(V)$. Since $\mathcal{O}(V)$ is integral over
$\mathcal{O}(V)^{\Gamma}$ there exists a polynomial of minimal
degree $\sum_{i=0}^n a_iX^i$ with $a_n=1$ and $a_i\in
\mathcal{O}(V)^{\Gamma}$ such that $\sum a_i u^i = 0$. Then
\[ 0 = \{ p, \sum_{i=0}^n a_iu^i \} = \sum_{i=0}^n (\{p,a_i\}u^i +
a_i \{ p,u^i\} ) = (\sum_{i=0}^n ia_iu^{i-1})\{ p, u\}.
\]
By minimality $\sum ia_iu^{i-1}\neq 0$, so $\{p,u\}= 0$. Thus
$p\in \text{Cas}(\mathcal{O}(V)) = \mathbb{C}$ as required.
\end{proof}

\subsection{Poisson structure} \label{sraPois} Henceforth we will concentrate
 on the case $t=0$. The results
 in this paragraph can be found in \cite{EG}. We let
 $A_{\bf{c}}$, respectively $Z_{\bf c}$, denote the algebra $A_{0,\bf{c}}$,
 respectively $Z_{0,\bf c}$, for any $\bf c\in \mathbb{A}^r$.
The family $Z_{\bf{c}}$ for $\bf{c} \in \mathbb{A}^r$, is flat,
and, by (\ref{projcong}), there is a $\C ^{\ast}$-action by
algebra automorphisms on the family lifting the natural action on
$\mathbb{A}^r$.

As we've already explained, in the degenerate case $c = {\bf 0}
\in \mathbb{A}^r$, the algebra $A_{\bf 0}\cong \mathcal{O}(V)\ast
\Gamma$ and $Z_{\bf 0} \cong \mathcal{O}(V)^{\Gamma}$, the skew
group ring and fixed point ring respectively. Thus, in view of
Proposition \ref{centre} and the discussion in the first paragraph
of (\ref{centre}) the variety $\text{Maxspec}(Z_{\bf c})$ is a
deformation of the quotient variety $V/\Gamma$.

Fix $\bf{c}\in \mathbb{A}^r$. The algebra $A_{\bf c}$ can be
lifted to a $\mathbb{C}$-algebra $\hat{A}_{\bf c}$ as described in
(\ref{alg1}) -- to be precise, define $\hat{A}_{\bf c}$ to be the
algebra $A_{t,{\bf c}}$, but with $t$ an indeterminate rather than
a complex number, so that $\hat{A}_{\bf c}$ is a $\C [t]$-algebra
with $\hat{A}_{\bf c}/t\hat{A}_{\bf c} \cong A_{\bf c}$. So by
(\ref{alg1}) $Z_{\bf c}$ admits a structure of Poisson algebra in
such a way that the pair $Z_{\bf c}\subseteq A_{\bf c}$ satisfies
the conditions of (\ref{setup}) -- that is, $A_{\bf c}$ is a
Poisson $Z_{\bf c}$-order. Recall the filtration (\ref{filt}) on
$A_{\bf c}$ with induced filtration on $Z_{\bf c}$, with
associated graded rings $\mathcal{O}(V) \ast \Gamma$ and
$\mathcal{O}(V)^{\Gamma}$ respectively.

There exists a $\C [t]$-algebra $\hat{Z}_{\bf c}$ such that the
flat family $Z_{\lambda \bf{c}}$ ($\lambda \in \C$) is realised by
specialisation, that is $Z_{\lambda {\bf c}} \cong \hat{Z}_{\bf c}
\otimes_{\C [t]} \C_{\lambda}$. There is, moreover, a Poisson
structure on $\hat{Z}_{\bf c}$ which is compatible with
specialisation and a $\C^{\ast}$-action on $\hat{Z}_{\bf c}$ (by
$\C$-algebra automorphisms) lifting the action on the $Z_{\lambda
{\bf c}}$'s.

\subsection{Symplectic leaves of $V/\Gamma$}\label{0case} The
restriction of the Poisson bracket on $\mathcal{O}(V)$ to
$\mathcal{O}(V)^{\Gamma}$ agrees with the Poisson bracket on the
Poisson $Z_{\bf 0}$-order described in the previous paragraph. We
shall determine the symplectic leaves of $V/\Gamma$. In
particular, we'll see that they are finite in number, so that
Proposition \ref{fteleaves} applies.

Given $v\in V$ let $\Gamma_v = \{ \gamma \in \Gamma :\gamma v =
v\}$, the stabiliser of $v$, and given $H\leq \Gamma$ let
$V^{o}_H= \{ v\in V:H = \Gamma_v \}$, and $V_H = \{ v\in V:H
\subseteq \Gamma_v \}$. Let $I(H) = \{ x^h-x : x\in
\mathcal{O}(V), h\in H\}$, an ideal of $\mathcal{O}(V),$ and set
\[ J(H) = I(H) \cap \mathcal{O}(V)^{\Gamma} = \bigcap_{\gamma \in \Gamma} I(H^{\gamma})
\cap \mathcal{O}(V)^{\Gamma},\] an ideal of
$\mathcal{O}(V)^{\Gamma}.$ Clearly $V_H$ is a closed subset of $V$
with $I(V_H) = I(H)$, and $V^o_H$ is open in $V_H,$ being the
complement of the closed subset of points with stabiliser strictly
containing $H$. Letting $H$ vary over subgroups of $\Gamma$ thus
gives a stratification of $V$ by locally closed subsets,
\[ V =\coprod_{H\leq \Gamma}V^{o}_H. \]
Let $\pi : V\longrightarrow V/\Gamma$ be the orbit map, and for
$H\leq \Gamma$ set $\mathcal{Z}^{o}_H = \pi(V^{o}_H)$, a locally
closed subset of $V/\Gamma$ which depends only on the conjugacy
class of $H$ in $\Gamma$. So there is a stratification of
$V/\Gamma$ by the locally closed sets $\mathcal{Z}^o_H$,
\begin{equation}\label{veegammstrat} V/\Gamma =\coprod_{H\leq
\Gamma}\mathcal{Z}^{o}_H,
\end{equation}
and
\[ \mathcal{Z}_H  := \overline{\mathcal{Z}}^o_H = \pi
(V_H),\] with $J(H)$ being the defining ideal of $\mathcal{Z}_H.$

\begin{prop} The symplectic leaves of $V/\Gamma$ are
precisely the sets $\mathcal{Z}^{o}_H$ as $H$ runs through the
conjugacy classes of subgroups of $\Gamma$ for which $V^o_H \neq
\emptyset.$
 The various leaves  $\mathcal{Z}^{o}_H$ coincide with the smooth points of the
irreducible components of the rank stratification of $V/\Gamma.$
\end{prop}
\begin{proof}
To show $J(H)$, and hence $\mathcal{Z}_H$, is Poisson closed take
$x,x'\in \mathcal{O}(V)$, $y\in \mathcal{O}(V)^\Gamma$ and $h\in H$.
Then, since the Poisson bracket is induced from the symplectic
form on $V$,
\begin{eqnarray*} \{ (x^h-x)x', y\} &=& \{x^h-x,y\}x' + (x^h-x)\{x',y\}
 \\&=& (\{x,y\}^h -\{x,y\})x'
+ (x^h-x)\{x',y\}, \end{eqnarray*} proving that $I(H)$, and
therefore $J(H)$, is stable under the Poisson action of
$\mathcal{O}(V)^\Gamma$. It follows from this and Proposition
\ref{alglf}.1 that the stratification of $V/\Gamma$ by symplectic
leaves is a refinement of (\ref{veegammstrat}).

Taking the union over subgroups conjugate to $H$ yields
\begin{equation} \label{cov}\pi^{-1}(\mathcal{Z}^{o}_H) = \coprod_{\gamma \in \Gamma} V^{o}_{H^{\gamma}}\end{equation}
and the restriction of $\pi$ to this space is a covering map whose
fibres have $[\Gamma :H]$ elements. Since $V^{o}_H$ is an open
subset of the vector space $V_H$, it follows that
\begin{equation*} \label{smooth} \mathcal{Z}^{o}_H \textit{ is
smooth. }
\end{equation*}
Moreover the restriction of the symplectic form to $V_H$ is
non-degenerate. Indeed, suppose $x\in V_H$ is in the radical. For
all $y\in {V}$, since $\Gamma \leq \text{Sp}(V)$,
$$(x,y) = |H|^{-1}\sum_{h\in H}(h.x,h.y) = (x, |H|^{-1}\sum_{h\in
H}h.y)= 0$$ ensuring $x=0$ as required. As a result, under the
covering (\ref{cov}), the form passes to the Poisson form on
$\mathcal{Z}^{o}_H$, proving non-degeneracy of the restriction of
the form to $\mathcal{Z}^o_H$. Thus each point of
$\mathcal{Z}^{o}_H$ has rank equal to $\dim (\mathcal{Z}^{o}_H)$
and so each subset $\mathcal{Z}^{o}_H$ is a symplectic leaf of
$V/\Gamma.$

In particular there are only finitely many leaves in $V/\Gamma$,
so that the second sentence of the proposition follows from
Proposition \ref{fteleaves}.2.
\end{proof}
\subsection{} It follows from the above the description and
Lemma \ref{D-M2} that the Poisson prime ideals of
$\mathcal{O}(V)^{\Gamma}$ are precisely the ideals $J(H)$, where
$H$ runs through representatives of conjugacy classes of
stabilisers of points in $V$. That this should be so was
speculated in \cite[Section 5]{AF}.
\subsection{Finite dimensional algebras}
\label{difffd} Given $v\in V$, recall that $\Gamma_v$ is its
stabiliser, and let $\mathfrak{m}_{\pi(v)}$ be the maximal ideal
of $\mathcal{O}(V)^{\Gamma}$ corresponding to $\pi(v)$. Recall the
finite dimensional algebra of (\ref{algIM}), $$(\mathcal{O}(V)\ast
\Gamma)_{\pi(v)} = \frac{\mathcal{O}(V)\ast
\Gamma}{\mathfrak{m}_{\pi(v)}\mathcal{O}(V)\ast \Gamma}.$$ Let $V=
V_{\Gamma_v} \oplus V'$ be a $\Gamma_v$-equivariant vector space
decomposition of $V$. It can be shown that $\Gamma_v$ is generated
by symplectic reflections on $V'$ and that there is an algebra
isomorphism
\begin{equation}
\label{gpim} (\mathcal{O}(V)\ast \Gamma)_{\pi(v)} \cong
\mat_{[\Gamma : \Gamma_v]} \left( (\mathcal{O}(V')\ast
\Gamma_v)_{\pi(0)}\right).
\end{equation}
In particular, this explicitly demonstrates that the
representation theory of $\mathcal{O}(V)\ast \Gamma$ is constant
along the symplectic leaves, as follows from Theorem \ref{algIM}.
\begin{rem} It can be shown quite generally that given an arbitrary
 vector space $V$ and finite group $\Gamma\leq GL(V)$, an exact analogue of
(\ref{gpim}) holds.
\end{rem}
\subsection{} \label{Alevqn}Let $W$ be a finite Weyl group. If $\mathfrak{h}$ denotes the reflection
representation of $W$, there is an induced action of $W$ by
symplectic reflections on $\mathfrak{h}\oplus \mathfrak{h}^*$. It
was suggested in \cite[Introduction]{AF} that the number of
Poisson prime ideals of $\mathcal{O}(\mathfrak{h}\oplus
\mathfrak{h}^*)^W$ with height $2k$ should equal $a_k$, the number
of conjugacy classes of $W$ having 1 as an eigenvalue with
multiplicity $k$ in the space $\mathfrak{h}$. \begin{prop} Retain
the above notation. The number of Poisson prime ideals in
$\mathcal{O}(\mathfrak{h}\oplus \mathfrak{h}^*)^W$ of height $2k$
equals the number of conjugacy classes of parabolic subgroups of
$W$ of rank $k$. In particular, this agrees with $a_k$ if and only
if all irreducible factors of $W$ are of type $A$.
\end{prop}
\begin{proof}
Without loss of generality we may assume that $W$ is irreducible,
and identify $\mathfrak{h}^*$ with $\mathfrak{h}$ via the Killing
form. Let $H$ be the stabiliser in $W$ of a point $(x,y)\in
\mathfrak{h}\oplus \mathfrak{h}$. It is well-known that $H$ is a
parabolic subgroup. Indeed, by definition $H$ is the intersection
the parabolic subgroups $W_x$ and $W_y$, \cite[Theorem
1.12]{humCOX}. To prove the claim we may assume that $x$ and $y$
are linearly independent in $\mathfrak{h}$. We will show that $H$
is the stabliser of $\lambda x + \mu y\in\mathfrak{h}$ for generic
values of $(\lambda : \mu ) \in \mathbb{P}^1$. Let $w\in W$, and
write $wx = \alpha x+\beta y + z$, where $x,y$ and $z$ are
linearly independent. Then $w$ stabilises $\lambda x + \mu
y\in\mathfrak{h}$ for a generic value of $(\lambda : \mu ) \in
\mathbb{P}^1$ if and only if
\[
wy = \frac{1}{\mu} ( \lambda(1-\alpha)x + (\mu- \lambda\beta)y -
\lambda z).
\]
By genericity we must have $\alpha =1$, $\beta =0$ and $z=0$,
proving that $w$ stabilises $x$ and hence $y$, as required.

Given a parabolic subgroup, $H$, the height of the corresponding
ideal $J(H)\subseteq \mathcal{O}(V)^{\Gamma}$ is the codimension
of the vector space $$(\mathfrak{h}\oplus\mathfrak{h}^*)_H = \{
(x,y): H \leq W_x \text{ and }H\leq W_y\} =\mathfrak{h}_H \oplus
\mathfrak{h}^*_H.$$ By \cite[1.15]{humCOX} $\dim\mathfrak{h}_H =
\dim \mathfrak{h} - \text{rank} (H)$, so the first claim of the
proposition follows.

Note that there is a well-defined injective mapping, say $\theta$,
from conjugacy classes of parabolic subgroups of $W$ to conjugacy
classes of $W$, sending a parabolic subgroup $H$ to the product of
its generating reflections, that is to a Coxeter element of $H$,
\cite[Proposition 3.1.15]{geckpfe}. Let $s$ be the rank function
on $W$, that is the function which assigns to an element $w$, the
codimension of $\{x \in \mathfrak{h} : w \in W_x\}$ in
$\mathfrak{h}$. In other words, $s(w)$ is the number of
eigenvalues of $w$ on $\mathfrak{h}$, not equal to 1. By
\cite[Theorem 1.12(d)]{humCOX}, $\dim \mathfrak{h}_H = \dim
\mathfrak{h} - s(\theta(H))$. Therefore, there is an one--to--one
association from Poisson prime ideals of
$\mathcal{O}(\mathfrak{h}\oplus \mathfrak{h}^*)^W$ to conjugacy
classes of elements in $W$, sending primes of height $2k$ to
elements of rank $k$.

The second claim of the proposition follows from \cite[Theorem
3.2.12, Section 3.4 and Appendix B]{geckpfe}: conjugacy classes of
parabolic subgroups provide a first approximation to conjugacy
classes of elements in $W$ via their Coxeter elements -- this
approximation is exact if and only if $W$ is of type $A$.
\end{proof}

\subsection{}\label{gencase}
We now show that for any $\bf{c}\in\mathbb{A}^r$ an analogue of
Proposition \ref{0case} holds for $\mathcal{Z}_{\bf c}
=\text{Maxspec} (Z_{\bf c})$.
\begin{trm}
For any $\bf{c}\in \mathbb{A}^r$ the symplectic leaves of
$\mathcal{Z}_{\bf c}$ are precisely the smooth points of the
irreducible components of the rank stratification. In particular
they are algebraic and finite in number.
\end{trm}
\begin{proof}
Let $\mathcal{X}\subseteq \mathcal{Z}_{\bf c}$ be a closed
subvariety and define a closed subvariety of
$\hat{\mathcal{Z}}_{\bf c}$ by $\hat{\mathcal{X}} =
\overline{\C^{\ast}\mathcal{X}}$, the Zariski closure of
$\C^{\ast}\mathcal{X}$. Taking the intersection of
$\hat{\mathcal{X}}$ with $\mathcal{Z}_0$ gives a closed subvariety
of $\mathcal{Z}_0$ which we denote $\text{gr}\mathcal{X}$. On the
level of ideals, this construction sends a radical ideal $I$ of
$Z_c$ to the ideal $\sqrt{\text{gr}I}$ of $Z_0$. Therefore $\dim
\mathcal{X} = \dim \text{gr}\mathcal{X}$ and if $\mathcal{X}$ is
Poisson closed, then so too is $\text{gr}\mathcal{X}$.

Let $r\in \mathbb{N}$ and suppose that some point
$\mathfrak{m}\in\mathcal{Z}_{\bf c}$ has rank $r$. Let
$\mathcal{U}$ be the subvariety of $\mathcal{Z}_{\bf c}$
consisting of the points whose rank is no more than $r$, and let
$\mathcal{X}$ be an irreducible component of $\mathcal{U}$. By
Lemma \ref{rank}.5, $\dim \mathcal{X} \geq r$. Now suppose that
$\mathfrak{n} \in \text{gr}\mathcal{X}$ has
$\text{rank}(\mathfrak{n}) = s > r$. Since $\hat{\mathcal{X}}$ is
Poisson closed, there is a non-empty open set of points of
$\hat{\mathcal{X}}$ whose rank is at least $s$. But, thanks to the
$\C^{\ast}$-action, $\hat{\mathcal{X}}$ is irreducible since
$\mathcal{X}$ is and has a non-empty open subset consisting of
points of rank $r$ since $\mathcal{X}$ does. This contradiction
shows that all points of $\text{gr}\mathcal{X}$ have rank at most
$r$. Therefore, since $\text{gr}\mathcal{X}$ is Poisson closed
$\dim \text{gr}\mathcal{X} \leq r$ by Proposition \ref{0case}.
Thus $$r\geq \dim \text{gr}\mathcal{X} = \dim \mathcal{X} \geq
r.$$ Thus the set $\mathcal{X}^o$ of points of $\mathcal{X}$ of
rank $r$ is non-empty, open, irreducible and of dimension $r$;
that is, it a symplectic manifold and hence is a leaf. The other
claims now follow immediately from Proposition \ref{fteleaves}.
\end{proof}
\subsection{}
\label{comm}  {\bf Remarks. 1.} Lemma \ref{gencase} shows that for
any ${\bf c}\in \mathbb{A}^r$ the pairs $Z_{\bf c}\subseteq A_{\bf
c}$ satisfy the hypotheses of Sections \ref{iss} and \ref{az}. In
general it is an interesting problem to describe the associated
families of finite dimensional algebras and the Azumaya
stratifications. Such a description for all $\bf{c}$ should among
other consequences give definitive information on
 the existence of a symplectic
resolution of $V/\Gamma$.

{\bf 2.} In the case $t \neq 0$ where $A_{t,{\bf c}}$ does {\em
not} satisfy a polynomial identity the most obvious open problem
is to describe the (two-sided) ideals of $A_{t,{\bf c}}$. The
mechanism of the proof of Proposition \ref{centre}.2 allows us to
induce the Poisson bracket on $\mathcal{O}(V)^{\Gamma}$ by
exploiting the fact that $\mathcal{O}(V) \ast \Gamma$ is the
associated graded algebra of $A_{t,{\bf c}}$. It therefore follows
easily that if I is an ideal of $A_{t,{\bf c}}$ then $\text{gr}I$
is a Poisson ideal of $A_{0,{\bf 0}}$, so that some information
about $I$ can be read off from Proposition \ref{0case}. But this
is unsatisfactory - one would prefer to associate $I$ with an
ideal of $A_{0,{\bf c}}$, and hence to make use of Proposition
\ref{gencase}. Whether this can be done, we leave as an open
question.

\end{document}